\documentclass[12pt]{amsart}
\usepackage{amssymb,amsmath, ulsy}
\usepackage[all]{xy}
\numberwithin{equation}{section}

\newcommand{\commentout}[1]{}

\newtheorem{thm}{Theorem}[section]
\newtheorem{lem}[thm]{Lemma}
\newtheorem{cor}[thm]{Corollary}

\newtheorem{ex}[thm]{Example}

\newcommand{\mc}[1]{\mathcal#1}

\newcommand{\Hom}{{\rm Hom}}

\newcommand{\SO}{\text{\rm SO}}     
\newcommand{\Sp}{\text{\rm Sp}}     
\renewcommand{\O}{\text{\rm O}}     
\newcommand{\U}{\text{\rm U}}       
\newcommand{\N}{\mathbb N}          
\newcommand{\R}{\mathbb R}          
\newcommand{\C}{\mathbb C}          
\newcommand{\Gr}{\text{\rm Gr}}     
\newcommand{\Rad}{\text{\rm Rad}}   
\newcommand{\proj}{\mathbf P}       
\newcommand{\mtx}[1]{\left[\,\begin{matrix} #1\end{matrix}\,\right ]}       
\newcommand{\smtx}[1]{\left[\begin{smallmatrix} #1\end{smallmatrix}\right ]}    
\renewcommand{\v}{\mathbf v}        
\renewcommand{\u}{\mathbf u}        
\newcommand{\w}{\mathbf w}          

\newcommand{\0}{\mathbf 0}          

\renewcommand{\i}{\mathbf i}
\newcommand{\F}{{\mathbb F}}
\newcommand{\Fc}{{\overline \F}}

\newcommand{\B}{\overline{B}}

\newcommand{\GL}{\text{\rm GL}}     

\begin{document}
\pagenumbering{arabic}
\title[Symmetric Subroup Actions on Isotropic Grassmannians]{Symmetric Subgroup Actions on Isotropic Grassmannians}
\author{Huajun HUANG and Hongyu HE}

\address{Department of Mathematics and Statistics, Auburn
University, Auburn, AL 36849} \email{huanghu@auburn.edu} 

\address{Department of Mathematics, Louisiana State University,
Baton Rouge, LA 70803} \email{hongyu@math.lsu.edu}

\footnote{This paper is partially supported by an NSF grant.}

\begin{abstract} Let $G$ be the group preserving a nondegenerate sesquilinear
form $B$ on a vector space $V$, and $H$ a symmetric subgroup of $G$
of the type $G_1\times G_2$. We explicitly parameterize the
$H$-orbits in $\Gr_G(r)$, the Grassmannian of $r$-dimensional
isotropic subspaces of $V$, by a complete set of $H$-invariants. We
describe the Bruhat order in terms of the majorization relationship
over a diagram of these $H$-invariants. The inclusion order, the
stabilizer, the orbit dimension, the open $H$-orbits, the
decompositions of an $H$ orbit into $H\cap G_0$ and $H_0$ orbits are
also explicitly described.
\end{abstract}
\maketitle 
\section{Introduction}
The symmetric subgroup orbits in flag manifolds have been
extensively studied. Their parametrization, in the most general
form, is due to Matsuki  \cite{Ma1, Ma2, Ma3, Ma4} and   Springer
\cite{Sp}.  There are finitely many such orbits.  In addition, there
is  a natural topological  ordering, called the Bruhat order,  among
the orbits: Orbits ${\mc O}\ge {\mc O}'$ if the Zariski closure of
${\mc O}$ contains ${\mc O}'$. The Bruhat order can be described
purely algebraically in terms of the Matsuki-Springer parameter
\cite{Ma2, rs, rs2, hel}.

In this paper, we thoroughly characterize one family of symmetric subgroup actions on the Grassmannians of isotropic subspaces by  complete sets of invariants, and we describe their Bruhat orders by majorization relationships over  diagrams of these invariants.
Consider the following symmetric pairs.

{\small
\begin{table}[h]
\caption{Symmetric pairs $(G,H)$ and  representation spaces $V$}
\vspace{-0.3 in}
\centering
$$
\begin{array}{|c|c|c|c|}
\hline
{\mathbf G}   &{\mathbf H}  &{\mathbf  V}
            &\text{\bf Conditions}
\\ \hline
\Sp_{2n}(\F)
    &\Sp_{2m}(\F)\times\Sp_{2n-2m}(\F)
        &\F^{2m}\oplus\F^{2n-2m}
            &0<m<n
\\ \hline
\O_n (\overline\F)
    &\O_m(\overline\F)\times\O_{n-m}(\overline\F)
        &\overline\F^m\oplus\overline\F^{n-m}
            &0<m<n
\\ \hline
\O(p,q)
    &\O(p_1,q_1)\times\O(p-p_1,q-q_1)
        &\R^{p_1+q_1}\oplus\R^{p-p_1+q-q_1}
            &0<p_1<p, \ 0<q_1<q
\\ \hline
\U(p,q)
    &\U(p_1,q_1)\times\U(p-p_1,q-q_1)
        &\C^{p_1+q_1}\oplus\C^{p-p_1+q-q_1}
            &0<p_1<p, \ 0<q_1<q
\\ \hline
\end{array}
$$
\begin{flushleft}
\text{ \ \ Note: $\F$ is an infinite field; $\overline\F$ is the algebraic closure of $\F$.}
\end{flushleft}
\label{symmetric-pairs}
\end{table}
}
\noindent
In Table \ref{symmetric-pairs}, $V$ is a  vector space  equipped with a {\em nondegenerate sesquilinear form} $B$; $G$ is the subgroup of $\GL(V)$ preserving the form $B$; $H$ is a symmetric subgroup of $G$ stabilizing two subspaces $U$ and $W$ of $V$, where $V=U\oplus W$ and $U\perp W$ with respect to form $B$.

A subspace $S$ of $V$ is said to be {\em isotropic} if $S\perp S$
with respect to form $B$. Let $\Gr_G(r)$ denote the projective
variety of $r$-dimensional isotropic subspaces of $V$. $\Gr_G(r)$ is
called an {\em isotropic Grassmannian} or simply {\em Grassmannian}.


This paper focuses on the $H$-action in $\Gr_G(r)$ for the triples
$(G,H,V)$ in Table \ref{symmetric-pairs}.

The symplectic case
$(G,H)=(\Sp_{2n}(\F),\Sp_{2m}(\F)\times\Sp_{2n-2m}(\F))$ has been
carefully investigated  \cite{RK, R}. Let $S$ be an isotropic
subspace of $V=U\oplus W$.
 In \cite{RK}, P. Rabau
and D. S. Kim construct an integral 4-tuple {\small
\begin{equation}\label{RK-invariant}
\left(\dim(S\cap U), \dim (S\cap W),
\dim\frac{S\cap(\proj_US\oplus\proj_W S)^\perp}{(S\cap
U)\oplus(S\cap W)},
\frac{1}{2}\dim\frac{S}{S\cap(\proj_US\oplus\proj_W S)^\perp}\right)
\end{equation}
}
to parameterize the $H$-orbit of $S$  and to determine the
  partial order induced by the inclusion of isotropic
subspaces. In \cite{R}, P. Rabau uses the 4-tuple to describe the
stabilizer of $S$ under the $H$-action and the codimension of an
$H$-orbit. In this paper, we remove the fraction $\frac 12$ from the
last term of \eqref{RK-invariant}  and denote the resulting 4-tuple
by  $(r_U(S),r_W(S),a(S),b(S))$. Then we determine the Bruhat order
in terms of the majorization relationship over a diagram of
$r_U(S)$, $r_W(S)$, $a(S)$ and $b(S)$.

In Section 2, we give a general discussion of $H\backslash\Gr_G(r)$
for the triples $(G,H,V)$
 in
Table \ref{symmetric-pairs}. The first main result is a
parametrization of the $H$-orbits  by a finite convex subset of an
integral lattice.

\begin{thm}\label{gen-orbit}
Let   $(G,H, V)$ be  given as in Table \ref{symmetric-pairs}.
\begin{enumerate}
  \item If $G=\Sp_{2n}(\F)$ or $\O_n(\overline{\F})$, then the
$H$-orbits in $\Gr_G(r)$
 can be parameterized by an integral 4-tuple of $H$-invariants $(r_U, r_W, a, b)$ defined by
 \eqref{H-inv}.
  \item If $G=\O(p,q)$ or $\U(p,q)$, then the $H$-orbits in $\Gr_G(r)$ can
be parameterized
 by an integral 5-tuple of $H$-invariants $(r_U, r_W, a, b_U, b_W)$, where $b_U$ and $b_W$ are defined in
\eqref{O(p,q)-b_U-b_W} and \eqref{U(p,q)-b_U-b_W}.
\end{enumerate}
In general, the $H$-orbit of   $S\in\Gr_G(r)$ is isomorphic to
$H/H_S$, where $H_S$ is the stabilizer of $S$ in the $H$-action. The
structure of $H_S$ is determined by  Theorem \ref{thm:stabilizer}.
The dimension of $H_S$ is given by Corollary
\ref{thm:orbit-stabilizer-dim}.
\end{thm}


The second main result we carry out is an explicit description of
the Bruhat order in terms of our $H$-invariants $(r_U, r_W, a, b)$
and $(r_U, r_W, a, b_U, b_W)$. The Bruhat order follows a simple
majorization relation over diagrams of these $H$-invariants.

\begin{thm}\label{gen-Bruhat}
 Consider the Bruhat order over $H\backslash \Gr_G(r)$ for $(G,H, V)$
 in Table \ref{symmetric-pairs}.
\begin{enumerate}
  \item When  $G=\Sp_{2n}(\F)$ or $\O_n(\overline{\F})$,
  we make the following diagram:
\begin{equation}\label{Bruhat-diagram-1}
\xymatrix@C=3pt @R=5pt {
    &b \ar@{-}[d]
\\
    &a \ar@{-}[dl] \ar@{-}[dr]
\\
r_U
    &   &r_W
}
\end{equation}
The  $H$-orbit parameterized by $(r_U,r_W,a,b)$
 is greater than the $H$-orbit parameterized by $(r_U',r_W',a',b')$
 in the Bruhat order   if and only if:
\begin{gather*}
\quad
 b\ge b',\quad a+b\ge a'+b',
 \quad
 r_U+a+b\ge r_U'+a'+b',\quad
r_W+a+b\ge r_W'+a'+b'.
\end{gather*}
Note that each of $b$, $a+b$, $r_U+a+b$ and $r_W+a+b$ is the sum of
all nodes connected to a given node via a descending path in diagram
\eqref{Bruhat-diagram-1}. The inequality $a+b\ge a'+b'$ is
redundant.
  \item When $G=\O(p,q)$ or $\U(p,q)$, we make the following
  diagram:
\begin{equation}\label{Bruhat-diagram-2}
\xymatrix@C=3pt @R=5pt {
b_U \ar@{-}[dr]
    &   &b_W \ar@{-}[dl]
\\
    &a \ar@{-}[dl] \ar@{-}[dr]
\\
r_U
    &   &r_W
}
\end{equation}
The  $H$-orbit parameterized by $(r_U,r_W,a,b_U,b_W)$
 is greater than the $H$-orbit parameterized by $(r_U',r_W',a',b_U',b_W')$
 in the Bruhat order   if and only if:
\begin{gather*}
b_U\ge b_U',\quad b_W\ge b_W',\quad a+b_U+b_W\ge a'+b_U'+b_W',
\\
 r_U+a+b_U+b_W\ge r_U'+a'+b_U'+b_W',\quad
 r_W+a+b_U+b_W\ge r_W'+a'+b_U'+b_W'.
\end{gather*}
Note that each of $b_U$, $b_W$, $a+b_U+b_W$, $r_U+a+b_U+b_W$ and
$r_W+a+b_U+b_W$ is the sum of all nodes connected to a given node
via a descending path in diagram \eqref{Bruhat-diagram-2}. The
inequality $a+b_U+b_W\ge a'+b_U'+b_W'$ is redundant.
\end{enumerate}
\end{thm}

Therefore, in terms of the Bruhat order, the $H$-invariants we use
provide the most natural way to describe the $H$-action in
$\Gr_G(r)$.

We also determine the inclusion order over the $H$-orbits of all
isotropic subspaces. Two $H$-orbits   ${\mc O}$ and ${\mc O}'$  have
the inclusion order ${\mc O}\succeq {\mc O}'$ if there exist $S\in
{\mc O}$ and $S'\in {\mc O}'$ such that $S\supseteq S'$. The third
main result of this paper is an extension of P. Rabau and D. S.
Kim's work on the inclusion order of symplectic case \cite[Theorem
4.3]{RK} to the inclusion orders of the other cases. See Theorems
\ref{thm:O(C)-inclusion-order}, \ref{thm:O(p,q)-inclusion-order} and
\ref{thm:U(p,q)-inclusion-order}.

When $G=\O(\C)$   or $G=\O(p,q)$, both $G$ and $H$ are disconnected.
Let $G_0$ be the identity component of $G$. We illustrate in
Sections 4 and 5 how an $H$-orbit in $\Gr_G(r)$ decomposes into
$(H\cap G_0)$ and $H_0$ orbits.

Our view point is purely algebraic. We derive our theorem by
analyzing the simultaneous isometry of a set of subspaces using the
tools presented in \cite[Theorem 5.3]{H}. It is unclear yet how our
parametrization should be identified with the Matsuki-Springer
parametrization \cite{bh, Ma1, Ma2, Sp}.

The $H$-orbits on isotropic Grassmannians  play an important
role in explicit construction of automorphic $L$-functions
\cite{gpr}.  The main motivation of this paper comes from representation theory.
Recall that functions on isotropic Grassmannian $\Gr_G(r)$ can be
used to define certain degenerate principal series $I_P(v)$.  The
representation $I_P(v)$ is one of the most intensively studied
series of representations.  In case $G/P$ is the Lagrangian
Grassmannian and $G$ the symplectic group, a preliminary
investigation by the second author gives a branching law for the
unitary $I_P(v)|_{H}$  \cite{com, ind}. This branching law is
multiplicity free and yields a Howe type $L^2$-correspondence
\cite{howe79, howe} between certain unitary representations of $G_1$ and
certain unitary representations of $G_2$ . So the remaining question is to
see if the degenerate principal series in other cases will decompose
in a similar fashion when restricted to $H$. A first step is to
understand how $H$ acts in $\Gr_G(r)$, in particular how $H$ acts on
the open orbits in $\Gr_G(r)$. The question of the structure of the
open orbits is answered in Corollaries \ref{symplectic:open},
\ref{O(C)-max-orbit}, \ref{O-max-orbit}, \ref{thm:U-max-orbit}.

Isotropic Grassmanian is a special case of partial flag variety.
Another interesting example of symmetric group action on flag
variety is the real semisimple group action on complex flag variety
\cite{Wo}. In this case, one often gets a finite number of open
orbits and the structure of these open orbits has broad implications
in representation theory.

\section{Preliminary}

\subsection{Settings}
Let
\begin{equation}
\label{N_0}
\N_0\ :=\ \{0\}\cup\N\ =\ \{0,1,2,3,\cdots\}.
\end{equation}
Let $\F$ be an infinite field, and $V$  a vector space over $\F$  equipped with a nondegenerate sesquilinear form $B$. Denote the {\em orthogonal direct sum}
$$
V=U\odot W
$$
if $V=U\oplus W$ as vector spaces and $U\perp W$ with respect to the form $B$, that is, $B(\u,\w)=0$ for any $\u\in U$ and $\w\in W$.
Suppose $\dim V=n$.
Let $G(V)$  or $G(n)$ denote the isometry  group of $V$ that preserves $B$:
$$
G(V)=G(n):=\{g\in\GL_{\F}(V)\mid B(g(\v),g(\v'))=B(\v,\v')\ \ \text{for}\ \ \v,\v'\in V\}.
$$
A symmetric pair  $(G,H)$ in Table \ref{symmetric-pairs} has the form
\ $(G(V), G(U)\times G(W))$ \ for certain sesquilinear form $B$ and certain decomposition
$V=U\odot W$.

For a fixed form $B$, let $\Gr_G(r)$ be the Grassmannian of $r$-dimensional isotropic subspaces of $V$.

%
%
%
%
%
%

\subsection{$H$-invariants in $\Gr_G(r)$}

Define the {\em radical} of a subspace $S$ of $V$ by
\begin{equation}
  \label{Rad}
  \Rad(S):= S\cap S^\perp=\{\v\mid \v\in S,\ B(\v,\v')=0 \text{ for any } \v'\in S\}.
\end{equation}
So $S$ is isotropic if and only if $\Rad(S)=S$.

A {\em flag} of a vector space is a nested sequence of subspaces.

Now suppose $S$ is an $r$-dimensional isotropic subspace of $V$,
namely $S\in\Gr_G(r)$. It induces a flag ${\mc F}_U(S)$ of $U$ and a
flag ${\mc F}_W(S)$ of $W$, respectively:
\begin{eqnarray}
\label{F_U(S)} \quad &{\mc F}_U(S): & \{\0\}\ \subseteq \ S\cap U\
\subseteq \ \Rad (\proj_U S)\ \subseteq \ \proj_U S
\\
\notag &&\qquad\qquad\qquad\qquad  \subseteq \ \Rad (\proj_U
S)^\perp\cap U\ \subseteq \ (S\cap U)^\perp \cap U \ \subseteq \ U.
\\
\label{F_W(S)} \quad &{\mc F}_W(S): & \{\0\}\ \subseteq \ S\cap W\
\subseteq \  \Rad (\proj_W S)\ \subseteq \ \proj_W S
\\
\notag &&\qquad\qquad\qquad\qquad
 \subseteq \ \Rad (\proj_W S)^\perp\cap W\ \subseteq \ (S\cap W)^\perp\cap W \ \subseteq \ W.
\end{eqnarray}
Here $\proj_U S:=(S+W)\cap U$ is the {\em projection} of $S$ onto
the $U$-component with respect to the decomposition $V=U\oplus W$.
Likewise for $\proj_W S$.

The symmetric subgroup $H=G(U)\times G(W)$ of $G$ consists of
elements of $\GL(V)$ that preserve  the form $B$ and stabilize the
subspaces $U$ and $W$. If  $h(S)=S'$ for $h\in H$ and
$S,S'\in\Gr_G(r)$, then $h$ sends each subspace in  ${\mc F}_U(S)$
(resp. ${\mc F}_W(S)$) bijectively to its counterpart in   ${\mc
F}_U(S')$ (resp.  ${\mc F}_W(S')$). So the dimensions of subspaces
in the flags ${\mc F}_U(\cdot)$ and ${\mc F}_W(\cdot)$ are
$H$-invariants.

For $\u\in\proj_U S$, let $\bar \u$ denote the element
$\u+\Rad(\proj_U S)$ in $\frac{\proj_U S}{\Rad(\proj_U S)}$. Then
$\frac{\proj_U S}{\Rad(\proj_U S)}$ has a nondegenerate sesquilinear
form
$\overline B$ induced by $B$:
\begin{equation}
  \label{induced-form}
  \overline B(\bar \u,\bar \u') := B(\u,\u'),\qquad \text{for \ } \bar\u,\ \bar\u'\in \frac{\proj_U S}{\Rad(\proj_U S)}.
\end{equation}
Similarly, $\frac{\proj_W S}{\Rad(\proj_W S)}$ has a nondegenerate sesquilinear form (also denoted by $\B$) induced by $B$.

Two vector spaces equipped with sesquilinear forms, $L_1$ with form
$B_1$ and $L_2$ with form $B_2$, are called {\em isometric}, if
there exists a linear bijection  $\phi:L_1\to L_2$, called an {\em
isometry}, such that $B_1( \u, \u')=B_2( \phi(\u), \phi(\u'))$ for
any $\u,\u'\in L_1$. 
With this notation, the isometry class of $\left(\frac{\proj_U S}{\Rad(\proj_U S)}, \B\right)$
is $H$-invariant   since $H$ preserves $B$. 


If $\u+\w,  \u'+\w'\in S$ such that $\u,\u'\in U$ and $\w,\w'\in W$, then
\begin{equation}
\label{dual-relation}
B(\u,\u')=-B(\w,\w').
\end{equation}
It immediately implies the following lemma.

\begin{lem}\label{thm:induced-form-correspondence}
The isometry class of
$\left(\frac{\proj_W S}{\Rad(\proj_W S)}, \B\right)$ is the additive inverse of the isometry class of $\left(\frac{\proj_U S}{\Rad(\proj_U S)}, \B\right)$.
\end{lem}


For example, if  $G\left(\frac{\proj_U S}{\Rad(\proj_U S)}\right )\simeq \O(p,q)$,
then $G\left(\frac{\proj_W S}{\Rad(\proj_W S)}\right)\simeq\O(q,p)$.

We are now ready to present  a complete set of $H$-invariants in
$\Gr_G(r)$ for symmetric pairs $(G,H)$ in Table
\ref{symmetric-pairs}. Denote
\begin{subequations}
\label{H-inv}
\begin{eqnarray}
\label{inv-r_U}
\qquad\
r_U(S) &:=& \dim(S \cap U),
\\
\label{inv-r_W}
r_W(S) &:=& \dim(S \cap W),
\\
\label{inv-a}
a(S) &:=&
\dim\frac{\Rad(\proj_U S)}{S \cap U}
= \dim\frac{\Rad(\proj_W S)}{S \cap W}
= \dim\frac{S\cap(\proj_US\odot\proj_W S)^\perp}{(S\cap U)\odot(S\cap W)},
\\
\label{inv-b}
b(S) &:=&
\dim\frac{\proj_U S}{\Rad(\proj_U S)}
=
\dim\frac{\proj_W S}{\Rad(\proj_W S)}
=
\dim\frac{S}{S\cap(\proj_US\odot\proj_W S)^\perp}.
\end{eqnarray}
\end{subequations}
Obviously, $r_U(S)$, $r_W(S)$, $a(S)$ and $b(S)$ are nonnegative integers and
\begin{equation}
\label{4-tuple-constraint}
r_U(S)+r_W(S)+a(S)+b(S)=\dim S=r.
\end{equation}

\begin{thm}\label{thm:invariant}
Let $(G,H,V)$ be a triple in Table \ref{symmetric-pairs}.
Let $r$ be an integer such that $0\le r\le \frac{1}{2}\dim V$.
For $S\in\Gr_G(r)$, the integral 4-tuple $(r_U(S), r_W(S), a(S), b(S))$ defined in \eqref{H-inv} and the isometry class of $\left(\frac{\proj_U S}{\Rad(\proj_U S)},\B\right)$
defined in \eqref{induced-form} form a complete set of $H$-invariants
that uniquely determines the $H$-orbit of $S$ in $\Gr_G(r)$.
\end{thm}

\begin{proof}
The $H$-invariant part is clear. Let us show that these
$H$-invariants uniquely determine the $H$-orbit of $S$ in
$\Gr_G(r)$. Let $S'\in\Gr_G(r)$ satisfy that
\begin{enumerate}
  \item $(r_U(S'), r_W(S'), a(S'), b(S')) = (r_U(S), r_W(S), a(S), b(S))$; 
  \item $\left(\frac{\proj_U S'}{\Rad(\proj_U S')}, \B\right)$ is isometric to
  $\left(\frac{\proj_U S}{\Rad(\proj_U S)}, \B\right)$. 
\end{enumerate}
We explicitly construct an  element of $H$ that sends $S$ to $S'$. For a subspace
$P_1$ of a vector space $P$, let
$$
P\ominus P_1
$$
denote the set of subspaces $P_2$ of $P$ such that $P=P_1\oplus P_2$.
\begin{enumerate}
  \item
  According to $r_U(S)=r_U(S')$, select a linear bijection \ $\phi_0:S\cap U\to S'\cap U$.

  \item
   Select $U_2\in \Rad(\proj_U S)\ominus (S\cap U)$ and $U_2'\in \Rad(\proj_U S')\ominus (S'\cap U)$.
  According to  $a(S)=a(S')$, select a linear bijection \ $\phi_1: U_2\to U_2'$.

  \item
  Select $U_3\in \proj_U S\ominus \Rad(\proj_U S)$. Then $(U_3, B)$ is isometric to
  $\left(\frac{\proj_U S}{\Rad(\proj_U S)},\B\right)$.
  Likewise, select $U_3'\in \proj_U S'\ominus \Rad(\proj_U S')$. Then $(U_3', B)$ is isometric to
  $\left(\frac{\proj_U S'}{\Rad(\proj_U S')},\B\right)$.
  Since $\left(\frac{\proj_U S'}{\Rad(\proj_U S')}, \B\right)$ and
  $\left(\frac{\proj_U S}{\Rad(\proj_U S)}, \B\right)$ are isometric,
  we can select an isometry \ $\phi_2: U_3\to U_3'$ \ with respect to the form $B$.

  \item
  Now according to \ $\proj_U S=(S\cap U)\odot U_2\odot U_3$, \
  the linear map $\phi:=\phi_0\oplus\phi_1\oplus\phi_2$ is an isometry from $\proj_U S$ to $\proj_U S'$. By Witt's extension theorem, $\phi$ can be extended to an isometry $h_U$ of $U$, that is, $h_U\in G(U)$.

  \item
  Next, according to $r_W(S)=r_W(S')$, select a linear bijection $\psi_0:S\cap W\to S'\cap W$.

  \item
  Select $W_{1,2}\in\proj_W S\ominus(S\cap W)$. Define a linear map $\widetilde{\psi_{1,2}}: W_{1,2}\to \proj_W S'$ as follow:
  Choose a basis $\{\w_1,\cdots,\w_k\}$ of $W_{1,2}$;
  for each $\w_i$, choose $\u_i\in\proj_U S$ such that $\u_i+\w_i\in S$;
  then $\phi(\u_i)\in\proj_U S'$;
  choose $\w_i'\in\proj_W S'$ such that $\phi(\u_i)+\w_i'\in S'$;
  define $\widetilde{\psi_{1,2}}(\w_i):=\w_i'$. It is routine to
  check that $\widetilde{\psi_{1,2}}$ is a well-defined linear injection
  and $W_{1,2}':=\text{Im}\; \widetilde{\psi_{1,2}}\in\proj_W S'\ominus(S'\cap W)$.
  Let $\psi_{1,2}:W_{1,2}\to W_{1,2}'$ be the linear bijection
  defined by $\psi_{1,2}(\w):=\widetilde{\psi_{1,2}}(\w)$. Then
  $\psi_{1,2}$ is an isometry  by \eqref{dual-relation}.

  \item
  Now according to $\proj_W(S)=(S\cap W)\odot W_{1,2}$, the linear map
  $\psi:=\psi_0\oplus\psi_{1,2}$ is an isometry from $\proj_W(S)$ to $\proj_W(S')$.
  By Witt's extension theorem, $\psi$ can be extended to an isometry
  $h_W\in G(W)$.

  \item
  Finally, $h:=h_U\times h_W$ is an element of $H$ that sends $S$ to $S'$.


\end{enumerate}

\end{proof}

\subsection{The stabilizer of $S\in\Gr_G(r)$ in the $H$-action}

Let $H_S$ denote the stabilizer of $S\in\Gr_G(r)$ in the $H$-action.
Then every $h\in H_S$ is of the form
$$
h=h_U\times h_W,
$$
where $h_U\in G(U)$ stabilizes the subspaces in the flag ${\mc
F}_U(S)$ defined in \eqref{F_U(S)}, and $h_W\in G(W)$ stabilizes the
subspaces in the flag ${\mc F}_W(S)$ defined in \eqref{F_W(S)}.

Choose a basis  ${\mc B}_U:=\{\hat\u_1,\cdots,\hat\u_{\dim U}\}$
 of $U$ such that:
\begin{enumerate}
  \item Each of the nontrivial subspaces in ${\mc F}_U(S)$, namely
$$
S\cap U,\quad
\Rad(\proj_U S),\quad
\proj_U S,\quad
\Rad(\proj_U S)^{\perp} \cap U,\quad
(S\cap U)^\perp\cap U,\quad\text{and \ }
U,
$$
is spanned by the first few vectors of ${\mc B}_U$.

\item
Note that $\frac{\proj_U S}{\Rad(\proj_U S)}$ and
$\frac{\Rad(\proj_U S)^{\perp} \cap U}{\Rad(\proj_U S)}$ are
nondegenerate with respect to their forms induced from $B$. We may
further assume that $\hat\u_i\perp \proj_U S$ for
$i=r_U+a+b+1,\cdots,\dim U-r_U-a$. Then
$$
\Rad(\proj_U S)^{\perp} \cap U
\ = \
\proj_U S\odot
\bigoplus_{i=r_U+a+b+1}^{\dim U-r_U-a}\F\hat\u_{i}.
$$
\end{enumerate}
With respect to the basis ${\mc B}_U$,
\begin{equation}
\label{h_U}
h_U=
\mtx{A_{11} &* &* &* &* &*
\\  &A_{22} &A_{23} &* &* &*
\\  & &A_{33} &0 &* &*
\\  & & &A_{44} &* &*
\\  & & & &A_{55} &*
\\  & & & & &A_{66}
}.
\end{equation}
Here
\begin{itemize}
  \item $A_{11}\in \GL_{r_U(S)}(\F)$ and $A_{66}\in \GL_{r_U(S)}(\F)$ uniquely determine each other;
  \item $A_{22}\in\GL_{a(S)}(\F)$ and $A_{55}\in\GL_{a(S)}(\F)$ uniquely determine each other;
  \item $A_{33}\in G\left(\frac{\proj_U S}{\Rad(\proj_U S)}\right)$;
  \item $A_{44}\in G\left(\frac{(\proj_U S)^{\perp}\cap U}{\Rad(\proj_U S)}\right)$.
\end{itemize}
Note that $h_U$ is in the parabolic subgroup $H(S,U)$  of $G(U)$
that preserves the flag
$$
\{\0\}\ \subseteq\  (S\cap U)\ \subseteq\  \Rad(\proj_U S)\
\subseteq\  \Rad(\proj_U S)^\perp\cap U\ \subseteq\ (S\cap
U)^\perp\cap U\ \subseteq\  U
$$
and has the Levi factor $\GL_{r_U(S)}(\F)\times\GL_{a(S)}(\F)\times
G\left(\frac{\Rad(\proj_U S)^{\perp} \cap U}{\Rad(\proj_U
S)}\right)$. The $G\left(\frac{\Rad(\proj_U S)^{\perp} \cap
U}{\Rad(\proj_U S)}\right)$ factor of $h_U$ in $H(S,U)$ is
$\mtx{A_{33} &0\\0 &A_{44}}.$

Next we consider $h_W\in G(W)$. In the basis ${\mc B}_U$ of $U$,
$$
(S\cap U)\odot\bigoplus_{i=1}^{a+b}\F\hat\u_{r_U+i}=\proj_U S.
$$
For each $\hat\u_{r_U+i}$ with $i=1,\cdots,a+b$, we select a vector in $\proj_W S$, denoted $\hat\w_{r_W+i}$, such that
$\hat\u_{r_U+i}+\hat\w_{r_W+i}\in S$. It is easy to see that
$$
(S\cap W)\odot\bigoplus_{i=1}^{a+b}\F\hat\w_{r_W+i}=\proj_W S.
$$
The set $\{\hat\w_{r_W+1},\cdots,\hat\w_{r_W+a+b}\}$ can be extended to a basis
${\mc B}_W:=\{\hat\w_1,\cdots,\hat\w_{\dim W}\}$ of $W$, such that:
\begin{enumerate}
  \item Each of the subspaces in ${\mc F}_W(S)$, namely
$$
S\cap W,\quad
\Rad(\proj_W S),\quad
\proj_W S,\quad
\Rad(\proj_W S)^\perp\cap W,\quad
(S\cap W)^\perp\cap W,\quad
\text{and \ }W,
$$
is spanned by the first few vectors of ${\mc B}_W$.
  \item $\hat\w_i\perp \proj_W S$ for $i=r_W+a+b+1,\cdots, \dim W-r_W-a$, so that
$$
\Rad(\proj_W S)^{\perp} \cap W
\ = \
\proj_W S\odot
\bigoplus_{i=r_W+a+b+1}^{\dim W-r_W-a}\F\hat\w_{i}.
$$
\end{enumerate}
Then with respect to the basis ${\mc B}_W$,
\begin{equation}
  \label{h_W}
h_W=  \mtx{B_{11} &* &* &* &* &*
\\
 &B_{22} &B_{23} &* &* &*
\\
 & &B_{33} &0 &* &*
\\
 & & &B_{44} &* &*
\\
 & & & &B_{55} &*
\\
 & & & & &B_{66} }.
\end{equation}
Here
\begin{itemize}
  \item
$B_{11}\in \GL_{r_W(S)}(\F)$ and $B_{66}\in\GL_{r_W(S)}(\F)$ uniquely determine each other;
  \item
$B_{22}\in\GL_{a(S)}(\F)$ and $B_{55}\in\GL_{a(S)}(\F)$ uniquely determine each other;
  \item
$B_{33}\in G\left(\frac{\proj_W S}{\Rad(\proj_W S)}\right)$;
  \item
$B_{44}\in G\left(\frac{(\proj_W S)^{\perp}\cap W}{\Rad(\proj_W S)}\right)$;
  \item
The conditions $\hat\u_{r_U+i}+\hat\w_{r_W+i}\in S$ for $i=1,\cdots,a+b$ and
$h_U\times h_W\in H_S$ imply that
  $\mtx{B_{22} &B_{23}\\ &B_{33}}=\mtx{A_{22} &A_{23}\\ &A_{33}}$.
\end{itemize}
The $h_W$ is in the parabolic subgroup $H(S,W)$ of $G(W)$ that
preserves the flag
$$
\{\0\}\ \subseteq\  (S\cap W)\ \subseteq\  \Rad(\proj_W S)\
\subseteq\  \Rad(\proj_W S)^\perp\cap W\ \subseteq\ (S\cap
W)^\perp\cap W\ \subseteq\  W
$$
and has the Levi factor \ $\GL_{r_W(S)}(\F)\times\GL_{a(S)}(\F)
\times G\left(\frac{[\Rad(\proj_W S)]^{\perp} \cap W}{\Rad(\proj_W
S)}\right)$. The $G\left(\frac{[\Rad(\proj_W S)]^{\perp} \cap
W}{\Rad(\proj_W S)}\right)$ factor of $h_W$ is $\mtx{B_{33} &0\\0
&B_{44}}$.

\begin{thm}\label{thm:stabilizer}
 Let $(G,H, V)$ be a triple in Table \ref{symmetric-pairs} and $S\in\Gr_G(r)$.
Then $h\in H_S$ if and only if $h=h_U\times h_W$, where $h_U\in G(U)$ and $h_W\in G(W)$ satisfy that:
\begin{enumerate}
\item
$h_U$ is of the form \eqref{h_U} with respect to the basis ${\mc B}_U$ of $U$.

\item
$h_W$ is of the form \eqref{h_W} with respect to the basis ${\mc B}_W$ of $W$.

\item
$h_U$ in \eqref{h_U} and $h_W$ in \eqref{h_W} are subjected to the constraint:
\begin{equation}
\label{h-U-W}
\mtx{A_{22} &A_{23}\\   &A_{33}}
\quad =\quad
\mtx{B_{22} &B_{23}\\   &B_{33}}.
\end{equation}
\end{enumerate}
\end{thm}

\begin{cor}
  \label{thm:orbit-stabilizer-dim} The dimension  of the $H$-orbit   ${\mc
  O}_S$ of $S$ is
\begin{equation}
 \dim {\mc O}_S = \dim G(U)+\dim G(W) - \dim H_S,
\end{equation}
and \ $\dim H_S$ \ equals to 
{\small
\begin{equation*}
\begin{split}
&\frac{1}{2} \left[\dim G(U)+\dim G(W) -\dim
G\left(\frac{\Rad(\proj_U S)^\perp\cap U}{\Rad(\proj_U S)}\right)-
  \dim G\left(\frac{\Rad(\proj_W S)^\perp\cap W}{\Rad(\proj_W
  S)}\right)\right]
\\
& \qquad
 +\dim G\left(\frac{(\proj_U S)^\perp\cap U}{\Rad(\proj_U
S)}\right)
  +\dim G\left(\frac{(\proj_W S)^\perp\cap W}{\Rad(\proj_W S)}\right)
+\dim G\left(\frac{\proj_U S}{\Rad(\proj_U S)}\right) 
\\
& \qquad \qquad +\frac{1}{2}\left[\dim \GL(S\cap U)+\dim\GL(S\cap
W)\right] -\dim\Hom_{\F}\left(\frac{\Rad(\proj_U S)}{S\cap
U},\frac{\proj_U S}{\Rad(\proj_U S)}\right).
\end{split}
\end{equation*}
}
\end{cor}

\begin{proof}
  It suffices to find the dimension of  Lie algebra of $H_S$. Let $\dim [A_{23}]$ denote the
  dimension of block $A_{23}$ when $h=h_U\times h_W$ goes through all
  elements of $H_S$. Then
\begin{gather*}
 \dim [A_{11}]=\dim\GL(S\cap U),\qquad \dim [A_{22}]=\dim
\GL\left(\frac{\Rad(\proj_U S)}{S\cap U}\right) , 
\\
\dim [A_{23}]=\dim\Hom_{\F}\left(\frac{\Rad(\proj_U S)}{S\cap
U},\frac{\proj_U S}{\Rad(\proj_U S)}\right), 
\\
 \dim[A_{33}]=\dim G\left(\frac{\proj_U
S}{\Rad(\proj_U S)}\right), \qquad 
\dim [A_{44}]=\dim G\left(\frac{(\proj_U S)^{\perp}\cap
U}{\Rad(\proj_U S)}\right).
\end{gather*}
 The other terms can be obtained easily. By the Levi decompositions of $H(S,U)$
and $H(S,W)$, and the constraints of $h=h_U\times h_W$ given by
Theorem \ref{thm:stabilizer}, it is straightforward to find $\dim
H_S$.
\end{proof}


\section{Symplectic Groups}\label{sect:symplectic}

Let $B$ be a nondegenerate symplectic form over $V\simeq\F^{2n}$. Let
$V=U\odot W$ where $\dim U=2m$ and $\dim W=2n-2m$. Then
\begin{equation}
\label{G-H-symplectic}
G=G(V)\simeq \Sp_{2n}(\F),\qquad H=G(U)\times G(W)\simeq \Sp_{2m}(\F)\times\Sp_{2n-2m}(\F).
\end{equation}
We first recall some results in \cite{RK, R} regarding the $H$-invariants and the stabilizer in  $\Gr_G(r)$ for integer $r$ with $0\le r\le n$. Then
 we shall study the Bruhat order of $H$-orbits in $\Gr_G(r)$.

\subsection{The $H$-invariants in $\Gr_G(r)$}

For $S\in\Gr_G(r)$, let $r_U(S), r_W(S), a(S)$ and $b(S)$ be defined in \eqref{H-inv}.
Theorem \ref{thm:invariant} verifies the following result of D. S. Kim and P. Rabau:

\begin{thm}\label{thm:symplectic}
\cite[Theorem 4.3]{RK}
 The $\Sp$-type $\left(r_U(S), r_W(S), a(S), \frac 12 b(S)\right)$
is a complete set of $H$-invariants that uniquely determines the
$H$-orbit of $S$ in $\Gr_G(r)$.
\end{thm}

Let $\left(r_U(S), r_W(S), a(S), b(S)\right)$  parameterize
 the $H$-orbit of $S$. Denote the $H$-orbit by
 ${\mc O}\left(r_U(S), r_W(S), a(S), b(S)\right)$.
 The range of $\left(r_U(S), r_W(S), a(S), b(S)\right)$ is as follow:

\begin{thm}
\label{thm:symplectic-inv-range}
\cite[Theorem 4.3]{RK}
 A 4-tuple $\left(r_U, r_W, a, b\right)\in{\N_0}^4$ parameterizes an $H$-orbit in $\Gr_G(r)$  if and only if $b$ is even, $r_U+r_W+a+b =r$, and
\begin{subequations}~\label{ineqsym}
\begin{eqnarray}
\label{Sp-U-ineq}
r_U+a+\frac{b}{2} &\leq& m,
\\
\label{Sp-W-ineq}
r_W+a+\frac{b}{2} &\leq& n-m.
\end{eqnarray}
\end{subequations}
\end{thm}

\subsection{The Bruhat order of $H\backslash\Gr_G(r)$}

The Bruhat order of the $H$-orbits in $\Gr_G(r)$ can be described by
elementary linear algebra method. The idea is that if ${\mc
O}\subseteq \Gr_G(r)$ and $S'$ is in the Zariski closure  $\overline
{\mc O}$, then for any
 subspace decomposition $V=R\oplus L$,
\begin{equation}
\label{Bruhat-idea} \overline{\lim_{S\in {\mc O}}}\; \dim \proj_R S
\ge \dim \proj_R S' , 
\qquad
 \overline{\lim_{S\in {\mc O}}}\;  \dim\frac{\proj_R S}{\Rad\left(\proj_R S\right)} \ge \dim
\frac{\proj_R S'}{\Rad\left(\proj_R S'\right)}.
\end{equation}

We make the
following diagram of $H$-invariants for $S\in\Gr_G(r)$:
\begin{equation}
\label{Sp-H-diagram}
\xymatrix@C=0pt @R=8pt
{
    &b(S) \ar@{-}[d]
\\
    &a(S) \ar@{-}[dl] \ar@{-}[dr]
\\
r_U(S)
    &   &r_W(S)
}
\end{equation}
The Bruhat order of $H\backslash \Gr_G(r)$ is characterized by a majorization relationship over
 diagram \eqref{Sp-H-diagram}. Define a partial order on  the diagram, where
 nodes $A\ge A'$ if and only if there exists a descending path from  $A$ to $A'$.
For each node $A$, we add the values of all nodes no less than  node
$A$. The resulting quantities are:
\begin{eqnarray*}
b(S) &=& \text{dim of a maximal nondegenerate subspace of $\proj_U S$}
\\
     &=&  \text{dim of a maximal nondegenerate subspace of $\proj_W S$,}
\\
a(S)+b(S) &=& \dim \frac{\proj_U S\odot \proj_W S}{S}
           =  \dim \frac{\proj_U S}{S\cap U}
           =  \dim \frac{\proj_W S}{S\cap W},
\\
r_U(S)+a(S)+b(S) &=& \dim \proj_U S,
\\
r_W(S)+a(S)+b(S) &=& \dim \proj_W S.
\end{eqnarray*}

\begin{thm}\label{thm:Sp-Bruhat-order}
The Bruhat order ${\mc O}(r_U,r_W,a,b)\ge {\mc O}(r_U',r_W',a',b')$
holds  in $H\backslash \Gr_G(r)$ if and only if the following
inequalities hold:
\begin{subequations}
\label{Sp-order-cond}
\begin{eqnarray}
\label{Sp-Bruhat-b}
b&\ge& b',
\\
\label{Sp-Bruhat-a+b}
a+b&\ge& a'+b',
\\
\label{Sp-Bruhat-r_U+a+b}
r_U+a+b&\ge& r_U'+a'+b',
\\
\label{Sp-Bruhat-r_W+a+b}
r_W+a+b&\ge& r_W'+a'+b'.
\end{eqnarray}
\end{subequations}
The inequality \eqref{Sp-Bruhat-a+b} is  implied by $r_U+r_W+a+b=r=r_U'+r_W'+a'+b'$, \eqref{Sp-Bruhat-r_U+a+b}
and \eqref{Sp-Bruhat-r_W+a+b}.
\end{thm}

\begin{proof}
By \eqref{Bruhat-idea},   inequalities   \eqref{Sp-order-cond} are
necessary  for ${\mc O}(r_U,r_W,a,b)\ge {\mc O}(r_U',r_W',a',b')$.
To show that they are sufficient, we will prove several claims
associate to some basic operations over the node values of  diagram
\eqref{Sp-H-diagram} that preserve the inequalities in
\eqref{Sp-order-cond}. These basic operations allow us to change
from $(r_U,r_W,a,b)$ to $(r_U',r_W',a',b')$ whenever  inequalities
\eqref{Sp-order-cond} hold.

Fix a   basis $\{\u_1,\cdots,\u_{2m}\}$ of $U$ and a basis $\{\w_1,\cdots,\w_{2n-2m}\}$
 of $W$ such that:
\begin{equation*}
\mtx{B(\u_i,\u_j)}_{2m\times 2m}=\smtx{0 &1\\-1 &0}^{\oplus m},
\quad
\mtx{B(\w_i,\w_j)}_{(2n-2m)\times(2n-2m)}={\smtx{0 &1\\-1 &0}}^{\oplus (n-m)}.
\end{equation*} 
Let ${\mc O}(r_U,r_W,a,b)$ be an $H$-orbit.
The following vectors  form a basis ${\mc B}_S$ of an element $S$ of ${\mc O}(r_U,r_W,a,b)$:
\begin{eqnarray*}
\u_{2i-1}+\w_{2i}\quad\text{and}\quad \u_{2i}+\w_{2i-1},
&\quad\text{for}\quad&
i=1,\cdots,b/2,
\\
\u_{b+2i}+\w_{b+2i},
&\quad\text{for}\quad&
i=1,\cdots, a,
\\
\u_{b+2a+2i},
&\quad\text{for}\quad&
i=1,\cdots,r_U,
\\
\w_{b+2a+2i},
&\quad\text{for}\quad&
i=1,\cdots,r_W.
\end{eqnarray*}

In the following arguments, we will define a subspace $S_x$ for $x\in\F$   spanned by
all vectors in  ${\mc B}_S$ but a few vectors being replaced. It will be easy to verify that:
\begin{enumerate}
  \item
  The entries of the basis vectors of $S_x$ given below are polynomials of $x$.
  \item
  $S_x\in {\mc O}(r_U,r_W,a,b)$ for every $x\in\F-\{0\}$.
\end{enumerate}
Then $S_0$ is in the Zariski closure of ${\mc O}(r_U,r_W,a,b)$ since $\F$ is an infinite field.
 In particular, $S_0\in\Gr_G(r)$.
If $S_0\in {\mc O}(r_U', r_W', a', b')$, then ${\mc O}(r_U,r_W,a,b)\ge {\mc O}(r_U', r_W', a', b')$
in the Bruhat order.
We assume that the related $H$-orbits in the following claims always exist.

\begin{enumerate}
\item
Claim: ${\mc O}(r_U,r_W,a,b)\ge {\mc O}(r_U,r_W,a+2,b-2)$.

Applying \eqref{ineqsym} to ${\mc O}(r_U,r_W,a+2,b-2)$,
$$
r_U+a+\frac{b}{2}\le m-1,
\qquad
r_W+a+\frac{b}{2}\le n-m-1.
$$
Let $S_x$ for $x\in\F$ be constructed as follow:
\begin{equation*}
\begin{array}{|c|c|}
\hline
    \text{\bf Vector in ${\mc B}_S$ being replaced}
        &\text{\bf Replaced by vector}
\\ \hline
    \u_b+\w_{b-1}
        &x(\u_b+\w_{b-1})+\u_{2m}+\w_{2n-2m}
\\ \hline
\end{array}
\end{equation*}
Then $S_x\in {\mc O}(r_U,r_W,a,b)$ for $x\in\F-\{0\}$ and $S_0\in {\mc O}(r_U,r_W,a+2,b-2)$.

\item
Claim: ${\mc O}(r_U,r_W,a,b)\ge {\mc O}(r_U+2,r_W,a,b-2)$.

Applying \eqref{Sp-U-ineq} to ${\mc O}(r_U+2,r_W,a,b-2)$,
$$
r_U+a+\frac{b}{2}\le m-1.
$$
Let $S_x$ for $x\in\F$ be constructed as follow:
\begin{equation*}
\begin{array}{|c|c|}
\hline
    \text{\bf Vectors in ${\mc B}_S$ being replaced}
        &\text{\bf Replaced by vectors}
\\ \hline
    \u_{b-1}+\w_{b}
        &\u_{b-1}+x\w_{b}
\\ \hline
    \u_{b}+\w_{b-1}
        &x^2\u_b+x\w_{b-1}+\u_{2m}
\\ \hline
\end{array}
\end{equation*}
Then $S_x\in {\mc O}(r_U,r_W,a,b)$ for $x\in\F-\{0\}$ and $S_0\in{\mc O}(r_U+2,r_W,a,b-2)$.

\item
Claim: ${\mc O}(r_U,r_W,a,b)\ge {\mc O}(r_U,r_W+2,a,b-2)$. The proof is similar.

\item
Claim: ${\mc O}(r_U,r_W,a,b)\ge {\mc O}(r_U+1,r_W+1,a,b-2)$.

Let $S_x$ for $x\in\F$ be constructed as follow:
\begin{equation*}
\begin{array}{|c|c|}
\hline
    \text{\bf Vectors in ${\mc B}_S$ being replaced}
        &\text{\bf Replaced by vectors}
\\ \hline
    \u_{b-1}+\w_{b}
        &x\u_{b-1}+\w_{b}
\\ \hline
    \u_b+\w_{b-1}
        &\u_b+x\w_{b-1}
\\ \hline
\end{array}
\end{equation*}
Then $S_x\in {\mc O}(r_U,r_W,a,b)$ for $x\in\F-\{0\}$ and $S_0\in {\mc O}(r_U+1,r_W+1,a,b-2)$.

\item
Claim: ${\mc O}(r_U,r_W,a,b)\ge {\mc O}(r_U+1,r_W,a+1,b-2)$.

Applying \eqref{Sp-U-ineq} to ${\mc O}(r_U+1,r_W,a+1,b-2)$,
$$
r_U+a+\frac{b}{2}\le m-1.
$$
Let $S_x$ for $x\in\F$ be constructed as follow:
\begin{equation*}
\begin{array}{|c|c|}
\hline
    \text{\bf Vector in ${\mc B}_S$ being replaced}
        &\text{\bf Replaced by vector}
\\ \hline
    \u_{b}+\w_{b-1}
        &x(\u_b+\w_{b-1})+\u_{2m}
\\ \hline
\end{array}
\end{equation*}
Then $S_x\in {\mc O}(r_U,r_W,a,b)$ for $x\in\F-\{0\}$ and
$S_0\in {\mc O}(r_U+1,r_W,a+1,b-2)$.

\item
Claim: ${\mc O}(r_U,r_W,a,b)\ge {\mc O}(r_U,r_W+1,a+1,b-2)$. The proof is similar.

\item
Claim: ${\mc O}(r_U,r_W,a,b)\ge {\mc O}(r_U+1,r_W,a-1,b)$.

Let $S_x$ for $x\in\F$ be constructed as follow:
\begin{equation*}
\begin{array}{|c|c|}
\hline
    \text{\bf Vector in ${\mc B}_S$ being replaced}
        &\text{\bf Replaced by vector}
\\ \hline
    \u_{b+2a}+\w_{b+2a}
        &\u_{b+2a}+x\w_{b+2a}
\\ \hline
\end{array}
\end{equation*}
Then $S_x\in {\mc O}(r_U,r_W,a,b)$ for $x\in\F-\{0\}$ and
$S_0\in {\mc O}(r_U+1,r_W,a-1,b)$.

\item
Claim: ${\mc O}(r_U,r_W,a,b)\ge {\mc O}(r_U,r_W+1,a-1,b)$. The proof is similar.

\end{enumerate}

By Theorem \ref{thm:symplectic-inv-range}, the set of $4$-tuples $(r_U,r_W,a,b)$ that parameterize
$H$-orbits in $\Gr_G(r)$ consists of the integer points in a convex set.
If $(r_U,r_W,a,b)$ and $(r_U',r_W',a',b')$ parameterize two $H$-orbits in $\Gr_G(r)$ and  they satisfy the inequalities in \eqref{Sp-order-cond}, we can find a sequence of
$4$-tuples 
\begin{multline*}
(r_U,r_W,a,b)=(r_U^{(0)},r_W^{(0)},a^{(0)},b^{(0)}),\quad
(r_U^{(1)},r_W^{(1)},a^{(1)},b^{(1)}),\quad
\cdots
\\
\cdots,\quad (r_U^{(d)},r_W^{(d)},a^{(d)},b^{(d)})
=(r_U',r_W',a',b'),
\end{multline*}
such that ${\mc O}(r_U^{(i-1)},r_W^{(i-1)},a^{(i-1)},b^{(i-1)})\ge {\mc O}(r_U^{(i)},r_W^{(i)},a^{(i)},b^{(i)})$ by one of the above claims for $i=1,\cdots,d$.
Then ${\mc O}(r_U,r_W,a,b)\ge {\mc O}(r_U',r_W',a',b')$
and the sufficient part is proved.
\end{proof}

\begin{cor}\label{symplectic:open}
Let $(G, H)$ be the symplectic symmetric pair in Table \ref{symmetric-pairs}.

\begin{enumerate}
\item
When $r< \min(2m, 2n-2m)$, the unique open $H$-orbit in
$\Gr_G(r)$ is
$$
\begin{cases}
{\mc O}(0,0,0,r) &\text{if $r$ is even;}
\\
{\mc O}(0,0,1,r-1) &\text{if $r$ is odd.}
\end{cases}
$$

\item
When  $\min(2m, 2n-2m)\le r\le n$, the unique open $H$-orbit in
$\Gr_G(r)$ is
$$
\begin{cases}
\mc O(0,r-2m, 0, 2 m) &\text{if $m\le n-m$;}
\\
\mc O(r-2n+2m, 0,0, 2n-2 m) &\text{if $m\ge n-m$.}
\end{cases}
$$

\end{enumerate}
\end{cor}


\begin{ex}
   Let $G=\Sp_{4m+8}(\F)$ and $H= \Sp_{2m}(\F) \times \Sp_{2m+8}(\F)$. We describe the Bruhat order of the $H$-orbits in the maximal isotropic Grassmannian $\Gr_G(2m+4)$.
Here $n=2m+4$ and $r=2m+4$. By Theorem \ref{thm:symplectic-inv-range},
$(r_U,r_W,a,b)\in{\N_0}^4$ parameterizes an $H$-orbit in $\Gr_G(2m+4)$ if and only if
the following constraints hold:
  \begin{equation*}
    \begin{cases}
      b \text{ \ is even;}
      \\
      r_U+r_W+a+b=2m+4;
      \\
      r_U+a+\frac{b}{2}\le m;
      \\
      r_W+a+\frac{b}{2}\le m+4.
    \end{cases}
    \quad\Longrightarrow\quad
    \begin{cases}
      b\in\{0,2,4,\cdots,2m\};
      \\
      a=0;
      \\
      r_U=m-\frac{b}{2};
      \\
      r_W=m+4-\frac{b}{2}.
    \end{cases}
  \end{equation*}
By Theorem \ref{thm:Sp-Bruhat-order}, the Bruhat order of $H\backslash \Gr_G(2m+4)$ is:
$$
{\mc O}(0,4,0,2m)>{\mc O}(1,5,0,2m-2)>{\mc O}(2,6,0,2m-4)>\cdots >{\mc O}(m,m+4,0,0).
$$
\end{ex}

\begin{ex}\label{ex:Sp-8-4-4-3}
Let $G=\Sp_8(\F)$ and $H=\Sp_4(\F)\times \Sp_4(\F)$. Then $n=4$ and $m=2$.
Consider the Bruhat order of $H$-orbits in $\Gr_G(3)$, where $r=3$.
By Theorem \ref{thm:symplectic-inv-range},
$(r_U,r_W,a,b)\in{\N_0}^4$ parameterizes an $H$-orbit in $\Gr_G(3)$ if and only if
$$
\text{$b$ is even;}\quad
r_U+r_W+a+b=3;\quad
r_U+a+\frac{b}{2}\le 2;\quad
r_W+a+\frac{b}{2}\le 2.
$$
So $b=0$ or $b=2$. There are 6 $H$-orbits in $\Gr_G(3)$ parameterized by:
$$
(r_U,r_W,a,b)\ \in\ \{(1,1,1,0),\
(1,2,0,0),\
(2,1,0,0),\
(1,0,0,2),\
(0,1,0,2),\
(0,0,1,2)\}.
$$
By Theorem \ref{thm:Sp-Bruhat-order}, the Bruhat order of $H\backslash\Gr_G(3)$ is given by the following diagram:
$$
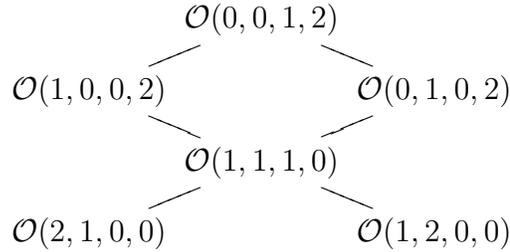
\xymatrix@C=0pt @R=8pt
{
    &{\mc O}(0,0,1,2) \ar@{-}[dl] \ar@{-}[dr]
\\
{\mc O}(1,0,0,2) \ar@{-}[dr]
    &   &{\mc O}(0,1,0,2) \ar@{-}[dl]
\\
    &{\mc O}(1,1,1,0) \ar@{-}[dl] \ar@{-}[dr]
\\
{\mc O}(2,1,0,0)
    &   &{\mc O}(1,2,0,0)
}
$$
\end{ex}

\subsection{The inclusion order of $H$-orbits}
In \cite{RK}, P. Rabau and D. S. Kim discuss an {\em inclusion
order} on the $H$-orbits of isotropic subspaces in all possible
dimensions, that is, on
$$
H\backslash\left(\bigcup_{r=0}^{n} \Gr_G(r)\right).
$$
The order is defined as follow:
$${\mc O}(r_U,r_W,a,b) \succeq {\mc O}(r_U',r_W',a',b')$$
if there exist $S\in {\mc O}(r_U,r_W,a,b)$ and $S'\in {\mc
O}(r_U',r_W',a',b')$ such that $S\supseteq S'$. Obviously, this
order is different from the Bruhat order, as any two distinct
$H$-orbits on a given $\Gr_G(r)$ have no ``$\succeq$'' relation.

\begin{thm}\label{thm:Sp-inclusion-order}
\cite[Theorem 4.3]{RK} For symplectic pair $(G,H)$ in Table
\ref{symmetric-pairs}, two $H$-orbits satisfy ${\mc O}(r_U,r_W,a,b)
\succeq {\mc O}(r_U',r_W',a',b')$ if and only if
$$
r_U\ge r_U',\quad
r_W\ge r_W',\quad
b\ge b',\quad
r_U+a+\frac{b}{2}\ge r_U'+a'+\frac{b}{2},\quad
r_W+a+\frac{b}{2}\ge r_W'+a'+\frac{b}{2}.
$$
\end{thm}

For $S\in {\mc O}(r_U,r_W,a,b)$,
$r_U=\dim (S\cap U)$, $b$ is the dimension of a maximal nondegenerate subspace of $\proj_U S$, and $r_U+a+\frac{b}{2}$ is the dimension of a maximal nilpotent subspace of $\proj_U S$.  Similarly for $r_W$ and $r_W+a+\frac{b}{2}$.

\subsection{Dimensions of orbit and stabilizer}

Let $S \in \mc O(r_U, r_W, a, b)\subseteq \Gr_G(r)$.  The stabilizer
$H_S$ of $S$ under the $H$-action is discussed in \cite[Section
5.II]{R}. The results are verified by Theorem \ref{thm:stabilizer}.

\begin{thm}
  \cite[Theorem 5.2]{R}
The codimension of $\mc O(r_U, r_W, a, b)$ is
{\small
\begin{eqnarray*}
&&\text{codim}\; \mc O(r_U, r_W, a, b)
\\
&=&
r_Ur_W+a(r_U+r_W)+2r_U(n-m-r_W-a-\frac{b}{2})+2r_W(m-r_U-a-\frac{b}{2})+\binom{a}{2}
\\
&=&
-3r_Ur_W-(r_U+r_W)(a+b)+2r_U(n-m)+2r_Wm+\binom{a}{2}.
\end{eqnarray*}
}
\end{thm}

By \ $ \dim \Gr_G(r)=2nr-\frac{3}{2}r^2+\frac{1}{2}r$, we get the
dimensions of $\mc O(r_U, r_W, a, b)$ and $H_S$:
\begin{eqnarray*}
\dim {\mc O}(r_U, r_W, a, b)
    &=&   \dim \Gr_G(r) - \text{codim}\; \mc O(r_U, r_W, a, b),
\\
\dim H_S
    &=& \dim H-\dim {\mc O}(r_U, r_W, a, b)
\\
    &=& \binom{2m+1}{2}+\binom{2n-2m+1}{2}-\dim {\mc O}(r_U, r_W, a, b).
\end{eqnarray*}
The results coincide with those of Corollary
\ref{thm:orbit-stabilizer-dim}.



\section{Orthogonal Groups on an Algebraically Closed Field}

Let $\Fc$ be an   algebraically closed field with $\text{char}(\Fc)\ne 2$.
Let $B$ be a nondegenerate symmetric form over $V\simeq\Fc^{n}$.
Suppose $V=U\odot W$
where $\dim U=m$ and $\dim W=n-m$. Then
\begin{equation}
  \label{G-H-orthogonal}
  G=G(V)\simeq\O_n(\Fc),
  \qquad
  H=G(U)\times G(W)\simeq\O_m(\Fc)\times\O_{n-m}(\Fc).
\end{equation}
We consider the $H$-action on the isotropic Grassmannian $\Gr_G(r)$  for
  $0\le r\le\lfloor n/2\rfloor$.

\subsection{The $H$-invariants in $\Gr_G(r)$}

Let  $S\in\Gr_G(r)$.
Then $\frac{\proj_U S}{\Rad(\proj_U S)}$ has a nondegenerate symmetric
   form $\B$ defined in \eqref{induced-form}.
The isometry class of $(\frac{\proj_U S}{\Rad(\proj_U S)},\B)$ is unique as $\Fc$ is algebraically closed.
Theorem \ref{thm:invariant}  implies the following result.

\begin{thm}\label{thm:orthogonal}
 The 4-tuple $\left(r_U(S), r_W(S), a(S), b(S)\right)$ defined in \eqref{H-inv}
is a complete set of $H$-invariants that determines the $H$-orbit of $S$ in $\Gr_G(r)$.
\end{thm}

Let $\left(r_U(S), r_W(S), a(S), b(S)\right)$  parameterize the
$H$-orbit of $S$ in $\Gr_G(r)$, and denote the $H$-orbit by ${\mc
O}\left(r_U(S), r_W(S), a(S), b(S)\right)$. The next theorem
determines the range of this 4-tuple. It is similar to Theorem
\ref{thm:symplectic-inv-range}. The proof is skipped.

\begin{thm}\label{thm:orthogonal-inv-range}
A 4-tuple $\left(r_U, r_W, a, b\right)\in{\N_0}^4$ parameterizes an $H$-orbit in $\Gr_G(r)$ if and only if:
\begin{enumerate}
  \item
  $r_U+r_W+a+b=r$.
  \item
  The following inequalities hold:
\begin{subequations}\label{O(C)-4-tuple-restriction}
\begin{eqnarray}
\label{O(C)-dim-U}
2r_U+2a+b
&\le& m,
\\ \label{O(C)-dim-W}
2r_W+2a+b
&\le& n-m.
\end{eqnarray}
\end{subequations}
\end{enumerate}
\end{thm}

%

Let $\{\u_1,\cdots,\u_m\}$  and
$\{\w_1,\cdots,\w_{n-m}\}$ be   orthonormal bases of $U$ and $W$, respectively.
Choose $\i\in \Fc$ such that $\i^2+1=0$.
The following vectors span a subspace of ${\mc O}(r_U,r_W,a,b)$
provided that $(r_U,r_W,a,b)$ satisfies the conditions in Theorem \ref{thm:orthogonal-inv-range}:
\begin{equation}
\label{O(C)-canonical-subspace}
\begin{cases}
\u_j+\i\w_j,
    &j=1,\cdots,b;
\\
\u_{b+2j-1}+\i\u_{b+2j}+\w_{b+2j-1}+\i\w_{b+2j},
    &j=1,\cdots,a;
\\
\u_{b+2a+2j-1}+\i\u_{b+2a+2j}
    &j=1,\cdots,r_U;
\\
\w_{b+2a+2j-1}+\i\w_{b+2a+2j}
    &j=1,\cdots,r_W.
\end{cases}
\end{equation}

\subsection{The Bruhat order of $H\backslash\Gr_G(r)$}

The Bruhat order of the $H$-orbits in $\Gr_G(r)$ is similar to that of the symplectic case. We make the following diagram of $H$-invariants:
\begin{equation}
\label{O(C)-H-diagram}
\xymatrix@C=0pt @R=8pt
{
    &b(S) \ar@{-}[d]
\\
    &a(S) \ar@{-}[dl] \ar@{-}[dr]
\\
r_U(S)
    &   &r_W(S)
}
\end{equation}
The  Bruhat order of $H\backslash \Gr_G(r)$ is characterized by the majorization relationship over
the diagram. For each node $A$ in the diagram, we define a quantity by adding the values of all nodes connected to node $A$ via descending paths. So we get $b(S)$, $a(S)+b(S)$, $r_U(S)+a(S)+b(S)$, and
$r_W(S)+a(S)+b(S)$.

\begin{thm}\label{thm:O(C)-Bruhat-order}
The Bruhat order ${\mc O}(r_U,r_W,a,b)\ge {\mc O}(r_U',r_W',a',b')$
holds in $\Gr_G(r)$ if and only if   the following inequalities
hold:
\begin{subequations}
\label{O(C)-order-cond}
\begin{eqnarray}
\label{O(C)-Bruhat-b}
b&\ge& b',
\\
\label{O(C)-Bruhat-a+b}
a+b&\ge& a'+b',
\\
\label{O(C)-Bruhat-r_U+a+b}
r_U+a+b&\ge& r_U'+a'+b',
\\
\label{O(C)-Bruhat-r_W+a+b}
r_W+a+b&\ge& r_W'+a'+b'.
\end{eqnarray}
\end{subequations}
The inequality \eqref{O(C)-Bruhat-a+b} is implied by \eqref{O(C)-Bruhat-r_U+a+b}
and \eqref{O(C)-Bruhat-r_W+a+b}.
\end{thm}

The proof is similar to that of Theorem \ref{thm:Sp-Bruhat-order} and we omit it here.

\begin{cor}\label{O(C)-max-orbit}
\text{}

\begin{enumerate}

\item
If $r\le\min\{m,n-m\}$,  the unique  open $H$-orbit  in $\Gr_G(r)$ is
${\mc O}(0,\ 0,\ 0,\ r)$.

\item
If $\min\{m,n-m\}\le r\le \lfloor n/2\rfloor$, the unique open $H$-orbit  in $\Gr_G(r)$  is
$$
\begin{cases}
{\mc O}(0,\ r-m,\ 0,\ m) &\text{if \ $\dim U\le \dim W$,}
\\
{\mc O}(r-n+m,\ 0,\ 0,\ n-m) &\text{if \ $\dim U> \dim W$.}
\end{cases}
$$
\end{enumerate}
\end{cor}

\begin{ex}
Let $G=\O_8(\Fc)$ and $H=\O_4(\Fc)\times\O_4(\Fc)$. Then $n=8$ and $m=4$.
Consider the Bruhat order of $H\backslash\Gr_G(3)$, where $r=3$.
By Theorem \ref{thm:orthogonal-inv-range}, $(r_U,r_W,a,b)\in{\N_0}^4$ parameterizes an $H$-orbit in $\Gr_G(3)$ if and only if
$$
r_U+r_W+a+b=3;\qquad 2r_U+2a+b\le 4;\qquad 2r_W+2a+b\le 4.
$$
Adding the last two inequalities and using the first equality, we have $a\le 1$. The possible 4-tuples
$(r_U,r_W,a,b)$ are:
$$
(2,1,0,0),\
(1,2,0,0),\
(1,1,0,1),\
(1,0,0,2),\
(0,1,0,2),\
(1,1,1,0),\
(0,0,1,2).
$$
By Theorem \ref{thm:O(C)-Bruhat-order}, the Bruhat order of $H\backslash\Gr_G(3)$ is given by the following diagram:
$$
\xymatrix@C=0pt @R=8pt
{
    &{\mc O}(0,0,1,2) \ar@{-}[dl] \ar@{-}[dr]
\\
{\mc O}(1,0,0,2) \ar@{-}[dr]
    &   &{\mc O}(0,1,0,2) \ar@{-}[dl]
\\
    &{\mc O}(1,1,0,1) \ar@{-}[d]
\\
    &{\mc O}(1,1,1,0) \ar@{-}[dl] \ar@{-}[dr]
\\
{\mc O}(2,1,0,0)
    &   &{\mc O}(1,2,0,0)
}
$$
Comparing with the case $G=\Sp_8(\F)$, $H=\Sp_4(\F)\times\Sp_4(\F)$ and $r=3$ in Example \ref{ex:Sp-8-4-4-3}, there is one additional orbit in the orthogonal case.
\end{ex}

\subsection{The inclusion order of $H$-orbits}
The inclusion  order  \ ${\mc O}(r_U,r_W,a,b) \succeq {\mc
O}(r_U',r_W',a',b')$ \ holds if there exist $S\in {\mc
O}(r_U,r_W,a,b)$ and $S'\in {\mc O}(r_U',r_W',a',b')$ such that
$S\supseteq S'$. This order for the case
$(G,H)=(\O_n(\Fc),\O_m(\Fc)\times\O_{n-m}(\Fc))$ is determined by
the same inequalities as in the symplectic case (cf. Theorem
\ref{thm:Sp-inclusion-order}).

\begin{thm}\label{thm:O(C)-inclusion-order}
The inclusion order of two $H$-orbits ${\mc O}(r_U,r_W,a,b) \succeq
{\mc O}(r_U',r_W',a',b')$ holds if and only if
\begin{equation}
\label{O(C)-inclusion-cond}
\begin{split}
r_U\ge r_U',\quad
r_W&\ge r_W',\quad
b\ge b',
\\
2r_U+2a+b\ &\ge\ 2r_U'+2a'+ b',
\\
2r_W+2a+ b\ &\ge\ 2r_W'+2a'+b'.
\end{split}
\end{equation}
\end{thm}

\begin{proof}
Suppose ${\mc O}(r_U,r_W,a,b) \succeq {\mc O}(r_U',r_W',a',b')$.
Let $S\in {\mc O}(r_U,r_W,a,b)$ and $S'\in¡¡{\mc O}(r_U',r_W',a',b')$
satisfy that $S\supseteq S'$. Then
\begin{itemize}
  \item $S\cap U\supseteq S'\cap U$;
  \item Every maximal nondegenerate subspace of $\proj_U S'$
  is contained in a maximal nondegenerate subspace of $\proj_U S$.
  \item Every minimal nondegenerate subspace of $U$  that contains
  $\proj_U S$   contains a minimal nondegenerate subspace of $U$
    that contains $\proj_U S'$;
  \item Similar arguments hold on the $W$ component.
\end{itemize}
By homogeneity property of  orbits, we may assume that $S$ is
spanned by the vectors in \eqref{O(C)-canonical-subspace}. Taking
dimensions, we get   inequalities  \eqref{O(C)-inclusion-cond}.

Conversely, suppose $(r_U,r_W,a,b)$ and $(r_U',r_W',a',b')$ satisfy
inequalities \eqref{O(C)-inclusion-cond}. Let $S\in {\mc
O}(r_U,r_W,a,b)$ be spanned by the vectors in
\eqref{O(C)-canonical-subspace}. If $a\ge a'$,   it is easy to find
a subspace $S'\subseteq S$ such that $S'\in {\mc
O}(r_U',r_W',a',b')$. Otherwise, $a<a'$. By a reduction process, we
may assume that $r_U'=0$, $r_W'=0$, $b'=0$ and $a=0$. Then
\eqref{O(C)-inclusion-cond} implies that \ $2\min\{r_U,r_W\}+b\ge
2a'$. Again, it is easy to find a subspace $S'$ of $S$ such that
$S'\in {\mc O}(r_U',r_W',a',b')= {\mc O}(0,0,a',0)$.
\end{proof}

\subsection{Dimensions of orbit and stabilizer}

Given  $S\in{\mc O}(r_U,r_W,a,b)$, Theorem \ref{thm:stabilizer} and
Corollary \ref{thm:orbit-stabilizer-dim} provide the structure of
$H_S$ as well as the dimensions of ${\mc O}(r_U,r_W,a,b)$ and $H_S$.

\begin{thm}\label{thm:O(C)-stabilizer}
For $S\in\mc O(r_U, r_W, a, b)$,
\begin{eqnarray*}
\dim {\mc O}(r_U, r_W, a, b)
=\dim H-\dim H_S
= \binom{m}{2}+\binom{n-m}{2}-\dim H_S,
\end{eqnarray*}
and \ $\dim H_S$ \ equals to 
{\small
\begin{equation} \label{dimor}
\begin{split}
 &\frac 12\left[\binom{m}{2}+ \binom{n-m}{2}+ r_U^2+r_W^2  - \binom{m-2r_U-2a}{2}- \binom{n-m-2r_W-2a}{2}\right]
\\
&\qquad +\binom{b}{2}-ab
+\binom{m-2r_U-2a-b}{2}+\binom{n-m-2r_W-2a-b}{2}.
\end{split}
\end{equation}
}
\end{thm}


\subsection{Decompose an $H$-orbit into $H\cap G_0$ and $H_0$ orbits}

When $\Fc=\C$, the identity components of $G$ and $H$ are, respectively:
$$
G_0=\SO_n(\Fc),
\qquad
H_0=\SO_m(\Fc)\times\SO_{n-m}(\Fc).
$$
Let $I_n$ be the $n\times n$ identity matrix. Denote the matrices
\begin{equation}
\label{I+-}
I_n^+:=I_n\qquad\text{and}\qquad
I_n^{-}:=I_{n-1}\oplus(-I_1).
\end{equation} 
The group $G$ has two connected components, namely
$G_0$ and $I_n^-G_0.$
For $t_1, t_2\in\{+,-\}$,
let $H^{t_1}_{t_2}$ denote
the connected component of $I_m^{t_1}\oplus I_{n-m}^{t_2}$ in $H$.
Then $H$ decomposes into two $(H\cap G_0)$-cosets and  four $H_0$-cosets as follow:
$$
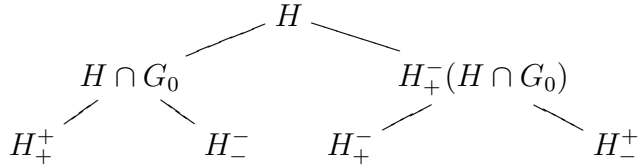
\xymatrix@C=2pt@R=6pt
{
    &    &   &H\ar@{-}[dll]\ar@{-}[drr]
                &   &   &
\\
    &H\cap G_0\ar@{-}[dl]\ar@{-}[dr]
         &   &  &   &H^{-}_{+}(H\cap G_0)\ar@{-}[dl]\ar@{-}[dr]
                        &
\\
H^+_+
    &   &H^-_-
            &   &H^-_+
                &   &H^+_-
                        &
}
$$

\begin{thm}\label{thm:O(C)-H-to-H-cap-G_0}
An $H$-orbit ${\mc O}(r_U,r_W,a,b)$ in $\Gr_G(r)$ always decomposes into 1 or 2 $(H\cap G_0)$-orbits.
Moreover, ${\mc O}(r_U,r_W,a,b)$  decomposes into 2 $(H\cap G_0)$-orbits if and only if both  equalities in \eqref{O(C)-dim-U} and \eqref{O(C)-dim-W} hold, if and only if
$r+a=n/2$.
\end{thm}

\begin{proof}
The number of $(H\cap G_0)$-orbits in the $H$-orbit of $S\in {\mc O}(r_U,r_W,a,b)$ equals to
$$[H:(H\cap G_0)H_S]\quad\in\quad\{1,2\}.$$
So $[H:(H\cap G_0)H_S]=2$ if and only if $(H\cap G_0)H_S=H\cap G_0$, if and only if
$H_S\subseteq  G_0$. Theorem \ref{thm:stabilizer} implies that the Levi factor of
$H_S$ is isomorphic to
$$
\GL_{r_U}\times\GL_{r_W}\times\GL_{a}\times G(b)\times G(m-2r_U-2a-b)\times G(n-m-2r_W-2a-b),
$$
where  for $K\in\{\GL_{r_U}, \GL_{r_W}, G(b)\}$, the $K$-component
of $H$ is diagonally embedded in a matrix group $K\times K$, and the
$\GL_a$-component is diagonally embedded in a matrix group
$\GL_a\times \GL_a\times \GL_a\times \GL_a$. So $H_S\subseteq
G_0=\SO_n(\Fc)$ if and only if
\begin{equation}\label{O(C)-4-tuple-equality}
m-2r_U-2a-b=0,\qquad
n-m-2r_W-2a-b=0.
\end{equation}
These are equivalent to the equalities in \eqref{O(C)-dim-U} and \eqref{O(C)-dim-W}.

It remains to prove the last statement.
If both equalities in \eqref{O(C)-dim-U} and \eqref{O(C)-dim-W} hold,
the sum of these two equalities produces $r+a=n/2$. Conversely, if $r+a=n/2$, then both equalities in \eqref{O(C)-dim-U} and \eqref{O(C)-dim-W} must hold.
\end{proof}

\begin{cor}
An open $H$-orbit of $\Gr_G(r)$ decomposes into 2 open $(H\cap G_0)$-orbits if and only if $r=n/2$,
in which $n$ is even and $\Gr_G(r)$ consists of maximal isotropic subspaces of $V$.
\end{cor}

\begin{thm}\label{thm:O(C)-H-cap-G_0-to-H_0}
The $(H\cap G_0)$-orbit of $S\in{\mc O}(r_U,r_W,a,b)$ decomposes into 1 or 2 $H_0$-orbits. Moreover,
it decomposes into 2 $H_0$-orbits if and only if $b=0$ and
at least one of the equalities in \eqref{O(C)-dim-U} and \eqref{O(C)-dim-W} holds.
\end{thm}

\begin{proof}
  The number of $H_0$-orbits in the $(H\cap G_0)$-orbit of $S$ equals to
$$
[H\cap G_0:H_0(H\cap G_0)_S]=[H\cap G_0:H_0(H_S\cap G_0)]\quad\in\quad\{1,2\}.
$$
Moreover, $[H\cap G_0:H_0(H_S\cap G_0)]=2$ if and only if $H_0(H_S\cap G_0)=H_0$,
if and only if $H_S\cap G_0\subseteq H_0$.
By the structure of $H_S$ described in Theorem \ref{thm:stabilizer},
$H_S\cap G_0\subseteq H_0$ if and only if  $b=0$, and at least one of the  equalities
in \eqref{O(C)-4-tuple-restriction} holds.
\end{proof}

\begin{cor}
An $H$-orbit ${\mc O}(r_U,r_W,a,b)$ in $\Gr_G(r)$ decomposes into 1, 2, or 4 $H_0$-orbits.
It decomposes into 4 $H_0$-orbits if and only if $b=0$ and $r+a=n/2$.
\end{cor}

\section{Real Orthogonal Groups}

Let $V:=\R^{p+q}$ be equipped with a
nondegenerate symmetric bilinear form $B$ of the type $I_{p}\oplus (-I_q)$.
Let $V=U\odot W$ be an orthogonal decomposition such that
$B|_U$ is of the type $I_{p_1}\times (-I_{q_1})$ and $B|_W$ is of the type $I_{p-p_1}\times (-I_{q-q_1})$. Denote
\begin{equation}
\label{O(p,q)-G-H}
G =G(V)\simeq\O(p,q),\qquad  H =G(U)\times G(W)\simeq \O(p_1,q_1)\times \O(p-p_1,q-q_1).
\end{equation}
We shall discuss the $H$-orbits in $\Gr_G(r)$ for
$0\le r\le\min\{p,q\}$.

\subsection{The $H$-invariants in $\Gr_G(r)$}

By Theorem \ref{thm:invariant}, the $H$-orbit of $S\in\Gr_G(r)$ is
determined by the $H$-invariants $r_U(S)$, $r_W(S)$, $a(S)$ and
$b(S)$ defined in \eqref{H-inv}
and by the  isometry class of
$(\frac{\proj_U S}{\Rad(\proj_U S)}, \B)$. According to
\eqref{dual-relation}, we  can denote
\begin{subequations}
\label{O(p,q)-b_U-b_W}
\begin{eqnarray}
\label{O(p,q)-b_U(S)}
\qquad
b_U(S)
&:=&
\text{the dimension of a maximal positive definite subspace of $\proj_U S$}
\\
\notag
&=&
\text{the dimension of a maximal negative definite subspace of $\proj_W S$,}
\\
\label{O(p,q)-b_W(S)}
b_W(S)
&:=&
\text{the dimension of a maximal negative definite subspace of $\proj_U S$}
\\ \notag
&=& \text{the dimension of a maximal positive definite subspace of $\proj_W S$.}
\end{eqnarray}
\end{subequations}
Then $b(S)=b_U(S)+b_W(S)$ and
$$
G\left(\frac{\proj_U S}{\Rad(\proj_U S)}\right)=\O\left(b_U(S),b_W(S)\right),
\quad
G\left(\frac{\proj_W S}{\Rad(\proj_W S)}\right)=\O\left(b_W(S),b_U(S)\right).
$$
The following result is  obvious.

\begin{thm}\label{thm:O(p,q)}
The 5-tuple $\left(r_U(S), r_W(S), a(S), b_U(S), b_W(S)\right)$
 is a complete set of $H$-invariants that determines the $H$-orbit of $S$ in $\Gr_G(r)$.
\end{thm}

We say that $\left(r_U(S), r_W(S), a(S), b_U(S), b_W(S)\right)$ parameterizes the $H$-orbit of $S$,
and the $H$-orbit  is denoted by
${\mc O}\left(r_U(S), r_W(S), a(S), b_U(S), b_W(S)\right)$.

\begin{thm}\label{thm:O(p,q)-inv-range}
A 5-tuple $(r_U,r_W,a,b_U,b_W)\in{\N_0}^5$ parameterizes an
$H$-orbit in $\Gr_G(r)$ if and only if the 5-tuple satisfies all
constraints below:
\begin{enumerate}
\item
$r_U+r_W+a+b_U+b_W=r$.

\item
The following conditions hold:
\begin{subequations}\label{O(p,q)-5-tuple-restriction}
\begin{eqnarray}
\label{O(p,q)-dim-nonneg-U}
r_U+a+b_U &\le & p_1,
\\ \label{O(p,q)-dim-nonpos-U}
r_U+a+b_W &\le & q_1,
\\ \label{O(p,q)-dim-nonneg-W}
r_W+a+b_W &\le & p-p_1,
\\ \label{O(p,q)-dim-nonpos-W}
r_W+a+b_U &\le & q-q_1.
\end{eqnarray}
\end{subequations}

\end{enumerate}
\end{thm}

\begin{proof}
First we prove the necessary part. Suppose that
$(r_U,r_W,a,b_U,b_W)$ parameterizes an $H$-orbit in $\Gr_G(r)$.
Obviously, $r_U+r_W+a+b_U+b_W=r$. Let $S^+$ be a maximal positive
definite subspace of $\proj_U S$. Then every vector $\v\in S_1:=
\Rad(\proj_U S)\odot S^+$ satisfies that $(\v,\v)\ge 0$. Let $U^-$
be a maximal negative definite subspace of $U$. Then $S_1\cap
U^-=\{\0\}$ and $S_1+ U^-\subseteq U$. So
$$
(r_U+a+b_U)+q_1=\dim S_1+\dim U^-=\dim (S_1\oplus U^-)\le \dim U=p_1+q_1.
$$
This leads to \eqref{O(p,q)-dim-nonneg-U}. Similarly for
\eqref{O(p,q)-dim-nonpos-U}, \eqref{O(p,q)-dim-nonneg-W}, and
\eqref{O(p,q)-dim-nonpos-W}.

Next we prove the sufficient part.
If a 5-tuple $(r_U,r_W,a,b_U,b_W)$ meets all the constraints in Theorem \ref{thm:O(p,q)-inv-range},
we find an isotropic subspace $S\in\Gr_G(r)$ whose $H$-orbit is parameterized by
$(r_U,r_W,a,b_U,b_W)$.
Let $\{\u_1^+,\cdots,\u_{p_1}^+,\u_1^-,\cdots, \u_{q_1}^-\}$
and
$\{\w_1^+,\cdots,\w_{p-p_1}^+,\w_1^-,\cdots,\w_{q-q_1}^-\}$ be  orthogonal bases of $U$ and $W$, respectively,
such that
\begin{equation}
\label{O(p,q)-basis} B(\u_i^+,\u_i^+)=1,\quad
B(\u_j^-,\u_j^-)=-1,\quad B(\w_t^+,\w_t^+)=1,\quad
B(\w_{\ell}^-,\w_{\ell}^-)=-1.
\end{equation}
Let $S$ be the subspace of $V$ spanned by the following basis
vectors:
\begin{gather}
\label{O(p,q)-S}
\left\{\u_i^+ +\w_i^-\right\}_{i=1}^{b_U}
\ \cup\
\left\{\u_i^- +\w_i^+\right\}_{i=1}^{b_W}
\ \cup\
\left\{\u_{b_U+i}^+ +\u_{b_W+i}^- +\w_{b_W+i}^+ +\w_{b_U+i}^-\right\}_{i=1}^a
\\ \notag
\ \cup\
\left\{\u_{b_U+a+i}^+ +\u_{b_W+a+i}^-\right\}_{i=1}^{r_U}
\ \cup\
\left\{\w_{b_W+a+i}^+ +\w_{b_U+a+i}^-\right\}_{i=1}^{r_W}.
\end{gather}
Then $S\in\Gr_G(r)$ and the $H$-orbit of $S$ is parameterized by
 $(r_U, r_W, a, b_U, b_W)$.
\end{proof}

\subsection{The Bruhat order of $H\backslash\Gr_G(r)$}

Construct the following diagram of $H$-invariants for $S\in\Gr_G(r)$:
\begin{equation}
\label{O(p,q)-H-diagram}
\xymatrix@C=0pt @R=8pt
{
b_U(S) \ar@{-}[dr]
    &   &b_W(S) \ar@{-}[dl]
\\
    &a(S) \ar@{-}[dl] \ar@{-}[dr]
\\
r_U(S)
    &   &r_W(S)
}
\end{equation}
The Bruhat order of $H\backslash \Gr_G(r)$ can be characterized by the majorization relationship over
this diagram. For each node $A$ in the  diagram, we define a quantity by adding the values of all nodes connected to node $A$ via descending paths. Then we get
\begin{equation*}
\begin{split}
&b_U(S)
=
\text{dim of a maximal positive definite subspace of $\proj_U S$,}
\\
&b_W(S)
=
\text{dim of a maximal negative definite subspace of $\proj_U S$,}
\\
&a(S)+b_U(S)+b_W(S)
=\dim\frac{\proj_U S\odot\proj_W S}{S}
=
\dim \frac{\proj_U S}{S\cap U}
=\dim \frac{\proj_W S}{S\cap W},
\\
&r_U(S)+a(S)+b_U(S)+b_W(S)
=
\dim\proj_U(S),
\\
&r_W(S)+a(S)+b_U(S)+b_W(S)
=
\dim\proj_W(S).
\end{split}
\end{equation*}

\begin{thm}\label{thm:O(p,q)-Bruhat-order}
The Bruhat order ${\mc O}(r_U,r_W,a,b_U,b_W)\ge {\mc
O}(r_U',r_W',a',b_U',b_W')$ holds in $\Gr_G(r)$ if and only if the
following inequalities hold:
\begin{subequations}\label{O(p,q)-order-cond}
\begin{eqnarray}
\label{O(p,q)-Bruhat-b_U}
b_U
&\ge&
b_U',
\\
\label{O(p,q)-Bruhat-b_W}
b_W
&\ge&
b_W',
\\
\label{O(p,q)-Bruhat-a+b_U+b_W}
a+b_U+b_W
&\ge&
a'+b_U'+b_W',
\\
\label{O(p,q)-Bruhat-r_U+a+b_U+b_W}
r_U+a+b_U+b_W
&\ge&
r_U'+a'+b_U'+b_W',
\\
\label{O(p,q)-Bruhat-r_W+a+b_U+b_W}
r_W+a+b_U+b_W
&\ge&
r_W'+a'+b_U'+b_W'.
\end{eqnarray}
\end{subequations}
Inequality \eqref{O(p,q)-Bruhat-a+b_U+b_W} is implied by
\eqref{O(p,q)-Bruhat-r_U+a+b_U+b_W} and
\eqref{O(p,q)-Bruhat-r_W+a+b_U+b_W}.
\end{thm}

\begin{proof}
The proof is similar to that of Theorem \ref{thm:Sp-Bruhat-order}.
The necessary part is obvious by \eqref{Bruhat-idea}. It remains to
prove the sufficient part. We will prove some claims associate to
several basic operations on the node values of  diagram
\eqref{O(p,q)-H-diagram} that preserve the order defined by
inequalities \eqref{O(p,q)-order-cond}. These operations allow us to
change from $(r_U,r_W,a,b_U,b_W)$ to $(r_U',r_W',a',b_U',b_W')$
whenever inequalities \eqref{O(p,q)-order-cond} hold.

Fix an orthogonal basis  $\{\u_1^+,\cdots,\u_{p_1}^+,\u_1^-,\cdots, \u_{q_1}^-\}$
  of $U$  and an orthogonal basis
$\{\w_1^+,\cdots,\w_{p-p_1}^+,\w_1^-,\cdots,\w_{q-q_1}^-\}$   of $W$ that satisfy \eqref{O(p,q)-basis}.
Let  $(r_U,r_W,a,b_U,b_W)$ be a  $5$-tuple that meets the constraints in Theorem \ref{thm:O(p,q)-inv-range}.

In the following arguments, we will define a subspace $S_x$ for $x\in\R$ spanned by all vectors in \eqref{O(p,q)-S} but one vector being replaced. It will be obvious that:
\begin{enumerate}
  \item
  The entries of the basis vectors of $S_x$ given below are polynomials of $x$.
  \item
  $S_x\in{\mc O}(r_U,r_W,a,b_U,b_W)$ for every $x\in\R^+$.
\end{enumerate}
Then $S_0$ is in the Zariski closure of ${\mc O}(r_U,r_W,a,b_U,b_W)$.
If $S_0\in {\mc O}(r_U',r_W',a',b_U',b_W')$, then ${\mc O}(r_U,r_W,a,b_U,b_W)\ge
{\mc O}(r_U',r_W',a',b_U',b_W')$ in the Bruhat order.

\begin{enumerate}

\item
Claim: ${\mc O}(r_U,r_W,a,b_U,b_W)\ge {\mc O}(r_U,r_W,a+1,b_U-1,b_W)$.

Applying \eqref{O(p,q)-dim-nonpos-U} and \eqref{O(p,q)-dim-nonneg-W} to $(r_U,r_W,a+1,b_U-1,b_W)$,
$$
r_U+a+b_W\le q_1-1,
\qquad
r_W+a+b_W\le p-p_1-1.
$$
Let $S_x$ be constructed as follow:
\begin{equation*}
\begin{array}{|c|c|}
\hline
    \text{\bf Vector in \eqref{O(p,q)-S} being replaced}
        &\text{\bf Replaced by the vector}
\\ \hline
    \u_{b_U}^++\w_{b_U}^-
        &(1+x)\u_{b_U}^++\u_{q_1}^-+\w_{p-p_1}^++(1+x)\w_{b_U}^-
\\ \hline
\end{array}
\end{equation*}
Then $S_x\in {\mc O}(r_U,r_W,a,b_U,b_W)$ for $x\in\R^+$ and $S_0\in {\mc O}(r_U,r_W,a+1,b_U-1,b_W)$.

\item
Claim: ${\mc O}(r_U,r_W,a,b_U,b_W)\ge {\mc O}(r_U+1,r_W,a,b_U-1,b_W)$.

Applying \eqref{O(p,q)-dim-nonpos-U}   to $(r_U+1,r_W,a,b_U-1,b_W)$,
$$
r_U+a+b_W\le q_1-1.
$$
Let $S_x$ be constructed as follow:
\begin{equation*}
\begin{array}{|c|c|}
\hline
    \text{\bf Vector in \eqref{O(p,q)-S} being replaced}
        &\text{\bf Replaced by the vector}
\\ \hline
    \u_{b_U}^++\w_{b_U}^-
        &(1+x^2)\u_{b_U}^++(1-x^2)\u_{q_1}^-+2x\w_{b_U}^-
\\ \hline
\end{array}
\end{equation*}
Then $S_x\in {\mc O}(r_U,r_W,a,b_U,b_W)$ for $x\in\R^+$ and $S_0\in {\mc O}(r_U+1,r_W,a,b_U-1,b_W)$.

\item
Claim: ${\mc O}(r_U,r_W,a,b_U,b_W)\ge {\mc O}(r_U,r_W+1,a,b_U-1,b_W)$. The proof is similar.

\begin{flushleft}
The next three claims are similar to the above three:
\end{flushleft}

\item
Claim: ${\mc O}(r_U,r_W,a,b_U,b_W)\ge {\mc O}(r_U,r_W,a+1,b_U,b_W-1)$.

\item
Claim: ${\mc O}(r_U,r_W,a,b_U,b_W)\ge {\mc O}(r_U+1,r_W,a,b_U,b_W-1)$.

\item
Claim: ${\mc O}(r_U,r_W,a,b_U,b_W)\ge {\mc O}(r_U,r_W+1,a,b_U,b_W-1)$.

\item
Claim: ${\mc O}(r_U,r_W,a,b_U,b_W)\ge {\mc O}(r_U+1,r_W,a-1,b_U,b_W)$.

Let $S_x$ be constructed as follows:
\begin{equation*}
\begin{array}{|c|c|}
\hline
    \text{\bf Vector in \eqref{O(p,q)-S} being replaced}
        &\text{\bf Replaced by the vector}
\\ \hline
    \u_{b_U+a}^++\u_{b_W+a}^-+\w_{b_W+a}^++\w_{b_U+a}^-
        &\u_{b_U+a}^++\u_{b_W+a}^-+x(\w_{b_W+a}^++\w_{b_U+a}^-)
\\ \hline
\end{array}
\end{equation*}
Then $S_x\in {\mc O}(r_U,r_W,a,b_U,b_W)$ for $x\in\R^+$ and $S_0\in {\mc O}(r_U+1,r_W,a-1,b_U,b_W)$.

\item
Claim: ${\mc O}(r_U,r_W,a,b_U,b_W)\ge {\mc O}(r_U,r_W+1,a-1,b_U,b_W)$. The proof is similar.

\end{enumerate}

These claims associate to some basic operations on the node values
of diagram \eqref{O(p,q)-H-diagram} that preserve the order defined
by   inequalities \eqref{O(p,q)-order-cond}. The last part is
similar to that of the proof of Theorem \ref{thm:Sp-Bruhat-order}.
\end{proof}

Theorem \ref{thm:O(p,q)-Bruhat-order} together with Corollary \ref{thm:O(p,q)-same-dim} implies the following result:

\begin{cor}\label{O-max-orbit}
\text{}

\begin{enumerate}
\item
When $r<\min\{p_1,q-q_1\}+\min\{q_1,p-p_1\}$, the open $H$-orbits in $\Gr_G(r)$ are  not unique,  and they are given by:
$${\mc O}(0,\ 0,\ 0,\ b_U,\ b_W),$$
where $b_U+b_W=r$, $0\le b_U\le\min\{p_1,q-q_1\}$ and  $0\le b_W\le\min\{q_1,p-p_1\}$.
\item
When $\min\{p_1,q-q_1\}+\min\{q_1,p-p_1\}\le r\le \min\{p,q\}$, there is a unique open $H$-orbit  in $\Gr_G(r)$ given by:
$$
\begin{cases}
{\mc O}(r-b_U-b_W,\ 0,\ 0,\ b_U,\ b_W) &\text{if \ $\dim U> \dim W$,}
\\
{\mc O}(0,\ r-b_U-b_W,\ 0,\ b_U,\ b_W) &\text{if \ $\dim U\le \dim W$,}
\end{cases}
$$
where $b_U=\min\{p_1,q-q_1\}$ and  $b_W=\min\{q_1,p-p_1\}$.
\end{enumerate}
\end{cor}

%

\begin{ex}
Let $G=\O(5,5)$ and $H=\O(2,3)\times\O(3,2)$. Then $p=5$, $q=5$,
$p_1=2$, $q_1=3$, $p-p_1=3$, $q-q_1=2$. Let $r=4$. We consider the
Bruhat order of $H$-orbits in $\Gr_G(4)$. By Theorem
\ref{thm:O(p,q)-inv-range}, an $H$-orbit ${\mc
O}(r_U,r_W,a,b_U,b_W)$ satisfies that:
\begin{gather*}
  r_U+r_W+a+b_U+b_W=4,\qquad
  r_U+a+b_U\le 2,\qquad
  r_U+a+b_W\le 3,
\\
  r_W+a+b_W\le 3,\qquad
  r_W+a+b_U\le 2.
\end{gather*}
Then $(r_U,r_W,a,b_U,b_W)$ could be one of the following 5-tuples:
\begin{gather*}
(0,0,0,1,3),\
(0,0,0,2,2),\
(0,1,0,1,2),\
(1,0,0,1,2),\
(1,1,0,0,2),
\\
(1,2,0,0,1),\
(2,1,0,0,1),\
(1,1,0,1,1),\
(2,2,0,0,0).
\end{gather*}
By Theorem \ref{thm:O(p,q)-Bruhat-order}, we obtain the following Bruhat order of $H\backslash \Gr_G(4)$:
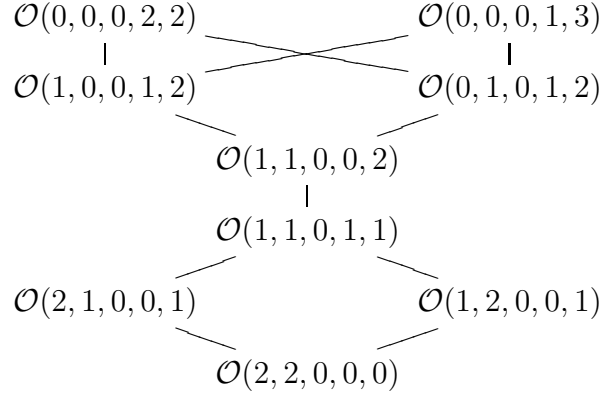
$$
\xymatrix@C=0pt @R=8pt
{
{\mc O}(0,0,0,2,2) \ar@{-}[d]  \ar@{-}[drr]
    &   &{\mc O}(0,0,0,1,3) \ar@{-}[dll]  \ar@{-}[d]
\\
{\mc O}(1,0,0,1,2) \ar@{-}[dr]
    &   &{\mc O}(0,1,0,1,2) \ar@{-}[dl]
\\
    &{\mc O}(1,1,0,0,2) \ar@{-}[d]
\\
    &{\mc O}(1,1,0,1,1) \ar@{-}[dl]  \ar@{-}[dr]
\\
{\mc O}(2,1,0,0,1) \ar@{-}[dr]
    &   &{\mc O}(1,2,0,0,1) \ar@{-}[dl]
\\
    &{\mc O}(2,2,0,0,0)
}
$$
There are two open $H$-orbits in this case.
\end{ex}

\subsection{The inclusion order of $H$-orbits}

The inclusion partial order ``$\succeq$'' for real orthogonal case is determined as follow:

\begin{thm}
  \label{thm:O(p,q)-inclusion-order}
There exist $S\in {\mc O}(r_U,r_W,a,b_U,b_W)$ and
$S'\in {\mc O}(r_U',r_W',a',b_U',b_W')$ such that $S\supseteq S'$ if and only if the following inequalities hold:
\begin{eqnarray*}
r_U\ge r_U',\qquad
r_W\ge r_W',
&\quad&
b_U\ge b_U',\qquad
b_W\ge b_W',
\\
r_U+a+b_U\ge r_U'+a'+b_U',
&\quad&
r_U+a+b_W\ge r_U'+a'+b_W',
\\
r_W+a+b_U\ge r_W'+a'+b_U',
&\quad&
r_W+a+b_W\ge r_W'+a'+b_W'.
\end{eqnarray*}
\end{thm}

The theorem can be proved by a similar reduction process as in the
proof of Theorem \ref{thm:O(C)-inclusion-order}.

\subsection{Dimensions of orbit and stabilizer}

Theorem \ref{thm:stabilizer} characterizes the structure of the
stabilizer $H_S$ of $S\in{\mc O}(r_U,r_W,a,b_U,b_W)\in\Gr_G(r)$.
Corollary \ref{thm:orbit-stabilizer-dim} implies the following
result.

\begin{thm}\label{thm:O(p,q)-stabilizer}
For $S\in O(r_U, r_W, a, b_U, b_W)$,
{\small
\begin{equation*}
\dim \mc O(r_U, r_W, a, b_U, b_W)
=\dim H-\dim H_S
= \binom{p_1+q_1}{2}+\binom{p-p_1+q-q_1}{2}-\dim H_S,
\end{equation*}
}
and \ $\dim H_S$ \ equals to
{\small
\begin{equation}
\label{O(p,q)-dim-H_S}
\begin{split}
&\quad\
\frac 12\binom{p_1+q_1}{2}
+\frac 12\binom{p+q-p_1-q_1}{2}
+\frac {r_U^2+r_W^2}2
+\binom{b_U+b_W}{2}
-a(b_U+b_W)
\\
&-\frac 12\binom{p_1+q_1-2r_U-2a}{2}
-\frac 12\binom{p+q-p_1-q_1-2r_W-2a}{2}
\\
&
+\binom{p_1+q_1-2r_U-2a-b_U-b_W}{2}
+\binom{p+q-p_1-q_1-2r_W-2a-b_U-b_W}{2}.
\end{split}
\end{equation}
}
\end{thm}

Formula \eqref{O(p,q)-dim-H_S} is similar to formula \eqref{dimor},
because $\dim_{\R}\O(p,q)=\dim_{\C}\O_{p+q}(\C)$ and
$\dim_{\R}\GL_a(\R)=\dim_{\C}\GL_a(\C).$ A direct consequence of
\eqref{O(p,q)-dim-H_S} is the following corollary.

\begin{cor}\label{thm:O(p,q)-same-dim}
If both $\mc O(r_U, r_W, a, b_U, b_W)$ and
$\mc O(r_U, r_W, a, b_U', b_W')$ exist and $b_U+b_W=b_U'+b_W'$, then
\ $\dim \mc O(r_U, r_W, a, b_U, b_W)=\dim \mc O(r_U, r_W, a, b_U', b_W').$
\end{cor}

\subsection{Decompose an $H$-orbit into $H\cap G_0$ and $H_0$ orbits}

We assume that $p, q, p_1, q_1, p-p_1, q-q_1>0$ for simplicity.
Denote the matrices $I_n^+:=I_{n}$ and $I_n^-:=I_{n-1}\oplus(-I_1)$.
Then $G=\O(p,q)$ has 4 connected components in Hausdorff topology, namely the $G_0$-cosets
of
$$
I_p^{t_1}\oplus I_q^{t_2}
\qquad\text{for}\qquad
t_1,\ t_2\in\{+,-\}.
$$
In particular, $\SO(p,q)=G_0\cup (I_p^-\oplus I_q^-)G_0$.
Similarly, $H=\O(p_1,q_1)\times\O(p-p_1,q-q_1)$ has 16 connected components, denoted by
$$
H^{t_1t_2}_{t_3t_4}\  := \
\left(I_{p_1}^{t_1}\oplus I_{q_1}^{t_2}\oplus I_{p-p_1}^{t_3}\oplus I_{q-q_1}^{t_4}\right)  H_0
\qquad\text{for}\quad t_1, t_2, t_3, t_4\in\{+,-\}.
$$
Obviously, $H^{++}_{++}=H_0$. Then $H$ and $H\cap G_0$ decompose into cosets as follow:
\begin{equation}
\label{O(p,q)-H-H-cap-G_0}
\xymatrix @C=0pt@R=1pc
{
    &   &   &   &   &H \ar@{-}[dlll] \ar@{-}[dl] \ar@{-}[dr] \ar@{-}[drrr]
\\
    &   &H\cap G_0 \ar@{-}[dll] \ar@{-}[dl] \ar@{-}[d] \ar@{-}[dr]
            &   &H^{--}_{++}(H\cap G_0) \ar@{.}[d]
                    &   &H^{+-}_{++}(H\cap G_0) \ar@{.}[d]
                        &   &H^{-+}_{++}(H\cap G_0) \ar@{.}[d]
\\
H_0
    &H^{--}_{--}
        &H^{+-}_{+-}
            &H^{-+}_{-+}
                &   &   &   &   &
}
\end{equation}

Fix an orthogonal basis
$$\{\u_1^+,\cdots,\u_{p_1}^+,\u_1^-,\cdots, \u_{q_1}^-\}\cup \{\w_1^+,\cdots,\w_{p-p_1}^+,\w_1^-,\cdots,\w_{q-q_1}^-\}$$
of $V$ that satisfies \eqref{O(p,q)-basis}.
Let $S\in {\mc O}(r_U,r_W,a,b_U,b_W)$ be the canonical subspace spanned by the vectors in \eqref{O(p,q)-S}.
The number of $(H\cap G_0)$-orbits in the $H$-orbit of $S$ equals to
$$
[H:(H\cap G_0)H_S]\in\{1,2,4\}.
$$
In fact, $[H:(H\cap G_0)H_S]=4/m$, where $m$ is the number of
$(H\cap G_0)$-cosets of $H$ that  intersect $H_S$. The number $m$
can be determined by
 $(r_U,r_W,a, b_U,b_W)$ as follows. Recall the constraints \eqref{O(p,q)-5-tuple-restriction}:
$$
\begin{array}{rclcrcl}
r_U+a+b_U
    &\le
        &p_1,
            &\quad
                &r_U+a+b_W
                    &\le
                        &q_1,
\\
r_W+a+b_W
    &\le
        &p-p_1,
            &\quad
                &r_W+a+b_U
                    &\le
                        &q-q_1.
\end{array}
$$

\begin{lem}\label{thm:O(p,q)-H_S-intersect-G-component}
The following statements hold:
\begin{enumerate}
\item
$H_S$ intersects  $H^{-+}_{++} (H\cap G_0)$ if and only if
$$r_U+a+b_U<p_1\qquad\text{or}\qquad r_W+a+b_W<p-p_1.$$

\item
$H_S$  intersects  $H^{+-}_{++}(H\cap G_0)$ if and only if
$$r_U+a+b_W<q_1\qquad\text{or}\qquad r_W+a+b_U<q-q_1.$$

\item
$H_S$ intersects \ $H^{--}_{++}(H\cap G_0)$ if and only if $a<r$.
\end{enumerate}
\end{lem}

\begin{proof}
The necessary part can be done by investigating the possible sign
combinations
 in the Levi factor of $h\in H_S$.

The sufficient part can be proved by the following explicit
construction:
\begin{enumerate}
\item
If  $r_W+a+b_W<p-p_1$,
 then $\w_{p-p_1}^+$ is not a component in the basis \eqref{O(p,q)-S} of $S$.
 So $\w_{p-p_1}^+\in S^\perp$ and
$$
I_{p_1}\oplus I_{q_1}\oplus I_{p-p_1}^-\oplus I_{q-q_1}
\quad \in\quad H_S\cap H^{++}_{-+}
\quad \subseteq\quad H_S\cap [H^{-+}_{++} (H\cap G_0)].
$$

Similarly,   $r_U+a+b_U<p_1$.

\item
The argument is similar to the preceding one.

\item
If $a<r$, then at least one of $b_U, b_W, r_U, r_W$ is greater than
0. Suppose $b_U>0$. Then $S$ has a basis vector $\u_1^+ +\w_1^-$ in
\eqref{O(p,q)-S}. Let $L\in\GL(V)$ have -1 eigenspace
$\text{span}\{\u_1^+,  \w_1^-\}$ and +1 eigenspace
$\text{span}\{\u_1^+,  \w_1^-\}^{\perp}$. Then
$$
L
\quad \in\quad
H_S\cap H^{-+}_{+-}
\quad \subseteq\quad
H_S\cap [H^{--}_{++}(H\cap G_0)]
$$ 
since $H^{-+}_{+-}=H^{--}_{++} H^{+-}_{+-}$. Similar for the other
cases.
\end{enumerate}
\end{proof}

Lemma \ref{thm:O(p,q)-H_S-intersect-G-component}
leads to the following result.

\begin{thm}\label{thm:O(p,q)-H-to-H-cap-G_0}
The number $N$ of  $(H\cap G_0)$-orbits in an $H$-orbit  is given
below: {\tiny
$$
\begin{array}{|c|c|c|c|p{3.6cm}|}
\hline
N   &G  &H  &\text{$H$-orbit}
                &\qquad\quad\text{Conditions}
\\ \hline
4   &\O(2a,2a)
        &\O(a,a)\times\O(a,a)
            &{\mc O}(0,0,a,0,0)
                &
\\ \hline
    &\O(p,2a)
        &\O(p_1,a)\times\O(p-p_1,a)
            &{\mc O}(0,0,a,0,0)
                &$p_1\ge a$, $p-p_1\ge a$, $p>2a$
\\ \cline{2-5}
    &\O(2a,q)
        &\O(a,q_1)\times\O(a,q-q_1)
            &{\mc O}(0,0,a,0,0)
                &$q_1\ge a$, $q-q_1\ge a$, $q>2a$
\\ \cline{2-5}
2   &
        &
            &
                &$r_U+a+b_U=p_1$,
\\
    &   &   &   &$r_W+a+b_W=p-p_1$,
\\
    &\O(p,q)
        &\O(p_1,q_1)\times\O(p-p_1,q-q_1)
            &{\mc O}(r_U,r_W,a,b_U,b_W)
                &$r_U+a+b_W=q_1$,
\\
    &(p=q)
        &   &   &$r_W+a+b_U=q-q_1$,
\\
    &   &   &   &$r_U+ r_W+ b_U+ b_W>0$
\\ \hline
1   &\multicolumn{4}{c|}{\text{all the other situations.}}
\\ \hline
\end{array}
$$
}
\end{thm}

Theorem \ref{thm:O(p,q)-H-to-H-cap-G_0}
and Corollary \ref{O-max-orbit}
imply the following decompositions of open $H$-orbits into open $(H\cap G_0)$-orbits.

\begin{cor}
An open $H$-orbit in $\Gr_G(r)$ always decomposes into 1 open
$(H\cap G_0)$-orbit except for the case $p=q=r$, in which the unique
open $H$-orbit
$$
\begin{cases}
 {\mc O}(0, r-p_1-q_1, 0, p_1, q_1)  &\text{if \ $p_1+q_1\le r$,}
\\
 {\mc O}(p_1+q_1-r, 0, r-q_1, r-p_1) &\text{if \ $p_1+q_1> r$,}
\end{cases}
$$
decomposes into 2 open $(H\cap G_0)$-orbits.
\end{cor}

\begin{proof}
Let ${\mc O}(r_U, r_W, a, b_U, b_W)$ be an open $H$-orbit in
$\Gr_G(r)$. By Corollary \ref{O-max-orbit}, we have $a=0$ and
$b_U+b_W>0$. According to Theorem \ref{thm:O(p,q)-H-to-H-cap-G_0},
${\mc O}(r_U, r_W, a, b_U, b_W)$ decomposes into 1 or 2 open $(H\cap
G_0)$-orbits, and it decomposes into   2   $(H\cap G_0)$-orbits if
and only if
\begin{gather*}
r_U+b_U=p_1,
\quad
r_W+b_W=p-p_1,
\quad
r_U+b_W=q_1,
\quad
r_W+b_U=q-q_1.
\end{gather*}
These together with $r= r_U+r_W+b_U+b_W$  imply that $p=q=r$. In
such case, Corollary \ref{O-max-orbit} gives the unique open
$H$-orbit.
\end{proof}


Similarly, the number of $H_0$-orbits in the $(H\cap G_0)$-orbit of $S\in\Gr_G(r)$ equals to
$$
[H\cap G_0: H_0(H\cap G_0)_S]=[H\cap G_0: H_0(H_S\cap G_0)]\quad\in\quad\{1,2,4\}.
$$
According to \eqref{O(p,q)-H-H-cap-G_0}, $[H\cap G_0:H_0(H_S\cap G_0)]=4/\ell$ where
$\ell$ is the number of cosets in $(H\cap G_0)/H_0=\{H_0, H^{--}_{--}, H^{+-}_{+-}, H^{-+}_{-+}\}$ that  intersect $H_S$. This can be determined by the 5-tuple $(r_U,r_W,a,b_U,b_W)$ and the sign combinations of the Levi factor of $h\in H_S$.
We omit the details here as there are many cases involved.

\section{Unitary Groups}

Suppose that $V:=\C^{p+q}$ is equipped with a nondegenerate
Hermitian form $B$ of the type $I_p\oplus(-I_q)$. Let $V=U\odot W$
be an orthogonal decomposition such that $B|_U$ is of the type
$I_{p_1}\times (-I_{q_1})$ and $B|_W$ is of the type
$I_{p-p_1}\times (-I_{q-q_1})$. Denote
\begin{equation}
\label{U(p,q)-G-H} G:=G(V)\simeq\U(p,q), \quad H:=G(U)\times
G(W)\simeq \U(p_1,q_1)\times\U(p-p_1,q-q_1).
\end{equation}
We consider the $H$-orbits in $\Gr_G(r)$ for $0\le r\le
\min\{p,q\}$. Most results in this section are similar to those in
Section 5. Their proofs are hence skipped. In counting the dimensions,
the subspaces of $V$ and the related quotient spaces   refer to
complex vector spaces, but all groups and orbits refer to real ones.

\subsection{The $H$-invariants in $\Gr_G(r)$}

Define
\begin{subequations}
\label{U(p,q)-b_U-b_W}
\begin{eqnarray}
\label{U(p,q)-b_U(S)} \qquad b_U(S) &:=& \text{the dimension of a
maximal positive definite subspace of $\proj_U S$}
\\
\notag &=& \text{the dimension of a maximal negative definite
subspace of $\proj_W S$,}
\\
\label{U(p,q)-b_W(S)} b_W(S) &:=& \text{the dimension of a maximal
negative definite subspace of $\proj_U S$}
\\ \notag
&=& \text{the dimension of a maximal positive definite subspace of
$\proj_W S$.}
\end{eqnarray}
\end{subequations}
Then  $b(S)=b_U(S)+b_W(S)$, and  the following result is true:

\begin{thm}\label{thm:U(p,q)-H-orbit}
The $5$-tuple $(r_U(S),r_W(S),a(S),b_U(S),b_W(S))$   is a complete
set of $H$-invariants that determines the $H$-orbit of $S$ in
$\Gr_G(r)$.
\end{thm}

Denote the $H$-orbit of $S\in\Gr_G(r)$ by ${\mc O}\left(r_U(S),
r_W(S), a(S), b_U(S), b_W(S)\right)$, where $\left(r_U(S), r_W(S),
a(S), b_U(S), b_W(S)\right)$  parameterizes the $H$-orbit.

\begin{thm}
\label{thm:U(p,q)-inv} A $5$-tuple $(r_U, r_W, a, b_U, b_W) \in
{\N_0}^5$ parameterizes an $H$-orbit in $\Gr_G(r)$ if and only if
the 5-tuple satisfies all constraints below:
\begin{enumerate}
\item
$r_U+r_W+a+b_U+b_W=r$.

\item
The following conditions hold:
\begin{subequations}\label{U(p,q)-5-tuple-restriction}
\begin{eqnarray}
\label{U(p,q)-dim-nonneg-U}
r_U+a+b_U &\le & p_1,
\\ \label{U(p,q)-dim-nonpos-U}
r_U+a+b_W &\le & q_1,
\\ \label{U(p,q)-dim-nonneg-W}
r_W+a+b_W &\le & p-p_1,
\\ \label{U(p,q)-dim-nonpos-W}
r_W+a+b_U &\le & q-q_1.
\end{eqnarray}
\end{subequations}
\end{enumerate}
\end{thm}

Let $G':=\O(p,q)$, $H':=\O(p_1,q_1)\times\O(p-p_1,q-q_1)$. Let $V'$
with a real symmetric form $B'$ be the natural representation space
of $(G',H')$. Let $\Gr_{G'}(r)$ be the $r$-dimensional isotropic
Grassmannian of $V'$. Apparently, there is a one-to-one
correspondence between $H\backslash \Gr_G(r)$ and $H'\backslash
\Gr_{G'}(r)$ for every $0\le r\le \min\{p,q\}$. So the Bruhat order
and the inclusion order here are the same as those in Section 5.


\subsection{The Bruhat order of $H\backslash\Gr_G(r)$}

Construct the following diagram of $H$-invariants for
$S\in\Gr_G(r)$:
\begin{equation}
\label{U(p,q)-H-diagram}
\xymatrix@C=0pt @R=8pt
{
b_U(S) \ar@{-}[dr]
    &   &b_W(S) \ar@{-}[dl]
\\
    &a(S) \ar@{-}[dl] \ar@{-}[dr]
\\
r_U(S)
    &   &r_W(S)
}
\end{equation}
The Bruhat order of $H\backslash \Gr_G(r)$ can be characterized by a
majorization relationship over this diagram, as in the real
orthogonal case.

\begin{thm}\label{thm:U(p,q)-Bruhat-order}
The Bruhat order ${\mc O}(r_U,r_W,a,b_U,b_W)\ge {\mc
O}(r_U',r_W',a',b_U',b_W')$ holds in $\Gr_G(r)$ if and only if  the
following inequalities hold:
\begin{subequations}\label{U(p,q)-order-cond}
\begin{eqnarray}
b_U
&\ge&
b_U',
\\
b_W
&\ge&
b_W',
\\
\label{U(p,q)-Bruhat-a+b_U+b_W} 
a+b_U+b_W &\ge& a'+b_U'+b_W',
\\
\label{U(p,q)-Bruhat-r_U+a+b_U+b_W} 
r_U+a+b_U+b_W &\ge&
r_U'+a'+b_U'+b_W',
\\
\label{U(p,q)-Bruhat-r_W+a+b_U+b_W}
r_W+a+b_U+b_W
&\ge&
r_W'+a'+b_U'+b_W'.
\end{eqnarray}
\end{subequations}
Moreover, inequality \eqref{U(p,q)-Bruhat-a+b_U+b_W} is implied by
inequalities \eqref{U(p,q)-Bruhat-r_U+a+b_U+b_W} and
\eqref{U(p,q)-Bruhat-r_W+a+b_U+b_W}.
\end{thm}

\begin{cor}\label{thm:U-max-orbit}
\text{}
\begin{enumerate}
\item
When $r<\min\{p_1,q-q_1\}+\min\{q_1,p-p_1\}$, the  open $H$-orbits in $\Gr_G(r)$ are not unique, and they are given by
$${\mc O}(0,\ 0,\ 0,\ b_U,\ b_W),$$
where $b_U+b_W=r$, $0\le b_U\le\min\{p_1,q-q_1\}$ and  $0\le b_W\le\min\{q_1,p-p_1\}$.
\item
When $\min\{p_1,q-q_1\}+\min\{q_1,p-p_1\}\le r\le \min\{p,q\}$, there is a unique open $H$-orbit in $\Gr_G(r)$ given by
$$
\begin{cases}
{\mc O}(r-b_U-b_W,\ 0,\ 0,\ b_U,\ b_W) &\text{if \ $\dim U> \dim W$,}
\\
{\mc O}(0,\ r-b_U-b_W,\ 0,\ b_U,\ b_W) &\text{if \ $\dim U\le \dim W$,}
\end{cases}
$$
where $b_U=\min\{p_1,q-q_1\}$ and  $b_W=\min\{q_1,p-p_1\}$.
\end{enumerate}
\end{cor}

\subsection{The inclusion order of $H$-orbits}

\begin{thm}
  \label{thm:U(p,q)-inclusion-order}
There exist $S\in {\mc O}(r_U,r_W,a,b_U,b_W)$ and
$S'\in {\mc O}(r_U',r_W',a',b_U',b_W')$ such that $S\supseteq S'$ if and only if the following inequalities hold:
\begin{eqnarray*}
r_U\ge r_U',\qquad
r_W\ge r_W',
&\quad&
b_U\ge b_U',\qquad
b_W\ge b_W',
\\
r_U+a+b_U\ge r_U'+a'+b_U',
&\quad&
r_U+a+b_W\ge r_U'+a'+b_W',
\\
r_W+a+b_U\ge r_W'+a'+b_U',
&\quad&
r_W+a+b_W\ge r_W'+a'+b_W'.
\end{eqnarray*}
\end{thm}

\subsection{Dimensions of orbit and stabilizer}

Theorem \ref{thm:stabilizer} gives the structure of ${\mc
O}(r_U,r_W,a,b_U,b_W)$. Corollary \ref{thm:orbit-stabilizer-dim}
together with the real dimensions
$$\dim \U(p,q)=(p+q)^2\quad\text{and}\quad\dim
\GL_{r_U}(\C)= 2r_U^2$$ 
gives the following explicit formulas of $\dim {\mc
O}(r_U,r_W,a,b_U,b_W)$ and $\dim H_S$ for any $S\in {\mc
O}(r_U,r_W,a,b_U,b_W)$.

\begin{thm}\label{thm:U(p,q)-stabilizer}
\text{}
\begin{eqnarray*}
\dim {\mc O}(r_U,r_W,a,b_U,b_W) = (p_1+q_1)^2+(p-p_1+q-q_1)^2-\dim
H_S
\end{eqnarray*}
and
\begin{equation}
\label{U(p,q)-dim-H_S}
\begin{split}
\dim H_S =
    &(p_1+q_1+r_W)^2+(p-p_1+q-q_1+r_U)^2-2(p+q)r+2r^2
\\
    &-4r_Ur_W+2a^2+(b_U+b_W)(2a+b_U+b_W),
\end{split}
\end{equation}
where \ $r=r_U+r_W+a+b_U+b_W$.
\end{thm}

\begin{cor}\label{U-same-dim}
If both $\mc O(r_U, r_W, a, b_U, b_W)$ and $\mc O(r_U, r_W, a, b_U',
b_W')$ exist and \ $b_U+b_W=b_U'+b_W'$, then \ $\dim \mc O(r_U, r_W,
a, b_U, b_W)=\dim \mc O(r_U, r_W, a, b_U', b_W')$.
\end{cor}

\vspace{3mm}
\begin{flushleft}
{\bf Acknowledgement:}  The authors would like to thank Professor R. Howe for his advice
on the manuscript, Professor G. Olafsson for pointing
out T. Matsuki's works, and the anonymous referee for
pointing to us the overlap with P. Rabau and D. S. Kim's work.
\end{flushleft}

\end{document}